\documentclass[12pt,leqno]{article}
\usepackage{amsmath,amssymb}
\usepackage[T1]{fontenc}
\usepackage{graphicx}
\usepackage{theorem}
\usepackage{manyfoot,perpage}
\DeclareNewFootnote{B}[fnsymbol]
\MakePerPage{footnoteB}

\theoremstyle{plain}
\newtheorem{theorem}{Theorem}

\DeclareMathOperator{\disc}{disc}

\newcommand{\fp}{\mathfrak p}
\newcommand{\Of}{\mathcal O}
\newcommand{\Z}{\ensuremath{\mathbb{Z}}}
\newcommand{\Q}{\ensuremath{\mathbb{Q}}}
\newcommand{\F}{\ensuremath{\mathbb{F}}}

\let\phi\varphi
\DeclareMathAlphabet{\mathbb}{U}{msb}{m}{n}

\title{Kurt Hensel on Common Inessential Discriminant Divisors, 1894}
\author{Fernando Q. Gouv\^ea and Jonathan Webster}
\date{}

\begin{document}


\vspace{\baselineskip}

\nocite{DirDed2}
\maketitle

The problem of the ``common inessential discriminant divisors''
attracted the attention of Dedekind, Kronecker, and Hensel in the
early days of algebraic number theory. Four sources are particularly
important: Dedekind's announcement, in 1871, of the second edition of
Dirichlet's lectures \cite{anzeige}, Dedekind's 1878 paper
\cite{Ded1878}, the 25th section of Kronecker's 1882
\textit{Grundz\"uge} \cite{Grundzuge}, and Hensel's 1894 paper
\cite{Hensel1894b}, which is our focus here. Both of the key papers of
Dedekind were translated and annotated in \cite{GW2}. A brief history
of results related to this problem can be found in
\cite[2.2.1]{Nark}. 

We here present an annotated translation of Kurt Hensel's
``Arithmetische Untersuchungen \"uber die gemeinsamen
ausserwesentlichen Discriminantentheiler einer Gattung''
(\textit{Journal f\"ur die Reine und Angewandte Mathematik},
\textbf{113} (1894), 128--160) \cite{Hensel1894b}. (All older volumes
of this journal are available online, so the original paper is easy to
access.) Our translation is based on a preliminary translation by
Timothy Molnar and Jonathan Webster, completed by Fernando
Q. Gouv\^ea, who also added explanatory footnotes. 

Before giving the translation itself we provide a quick outline of the
mathematical background of this paper. Some of what we say is informed
speculation: it seems clear that both Dedekind and Kronecker started
thinking about this subject early on, certainly by the 1860s, but
neither one published anything for quite a while: Dedekind published
his first version of the theory in 1871 \cite{SuppX} and Kronecker
finally explained his theory in the \textit{Grundz\"uge} of 1882
\cite{Grundzuge}. Few notes or unpublished manuscripts seem to have
survived. Any account of their process is inferred from what they said
later. For the evolution of Dedekind's ideas, see also
\cite{EdwardsIdeals}, \cite{Haubrich}, and \cite{Haffner}.

\section{The Mathematical Background}

When Kronecker and Dedekind set out to generalize Kummer's theory of
cyclotomic integers, they quickly ran into obstacles. Finding a way
around these difficulties led each of them to develop a far more
complicated theory than Kummer's. As a result, each had to justify the
extra work by highlighting what made it necessary.

Suppose $n>0$ is an integer and let $\zeta$ be a primitive $n$-th root
of unity. Kummer had found an explicit description in terms of
congruences of how rational primes factor in the cyclotomic integers
$\Z[\zeta]$. It seems that both Dedekind and Kronecker\footnote{And
  also Selling in \cite{Selling}.} saw that Kummer's description could
be interpreted in terms of congruences between polynomials (known as
``higher congruences'' at the time). In modern terms, it would go
something like this.

\begin{theorem} Let $n>2$ be an integer, let $\zeta$ be a
  primitive $n$-th root of unity, and let $\Phi_n(x)$ be the $n$-th
  cyclotomic polynomial. Fix a prime number $p\in\Z$ and let
  \[ \Phi_n(x)\equiv F_1(x)^{e_1} F_2(x)^{e_2} \dots F_r(x)^{e_r}
    \pmod p\] be the factorization of $\Phi(x)$ modulo $p$, where the
  $F_i(x)$ are distinct irreducible polynomials in $\F_p[x]$. Then the
  factorization of $(p)$ in $\Z[\zeta]$ is
  \[ (p)=\fp_1^{e_1}\fp_2^{e_2}\dots\fp_r^{e_r},\]
  with distinct prime ideals $\fp_i=(p,F_i(\zeta))$.
\end{theorem}

Of course, Kummer did not speak of ideals; instead, he thought of
$\fp_i$ as the ``ideal prime divisor'' determined by $p$ and
$F_i(x)$. He gave an explicit method for determining the exponent of
$\fp_i$ in a factorization. Thus, the ``ideal prime divisor'' is
essentially the valuation corresponding to $\fp_i$. 

This beautiful result seemed to suggest the possibility of a very
simple theory in the general case: for a general number field
$\Q(\alpha)$, let $\Phi(x)$ be the minimal polynomial for $\alpha$ and
factor it modulo $p$. One could then use this to define ``ideal
primes'' \`a la Kummer.

The choice of $\alpha$ is crucial, of course. At least one example
would have been familiar to everyone: the field $\Q(\sqrt{-3})$ is the
same as the cyclotomic field of order $3$. Kummer's approach worked if
one took $\alpha$ to be a cube root of $1$ but would not work if we
took $\alpha=\sqrt{-3}$. Both Dedekind and Kronecker figured out that
one needed to work with all the algebraic integers in the field
$\Q(\alpha)$.

That highlights the first difficulty: in the case of $\Q(\zeta)$ the
ring of algebraic integers is exactly $\Z[\zeta]$, but this will not
be true in general. If $K$ is a number field and $\Of\subset K$ is its
ring of algebraic integers there may not exist any $\alpha\in\Of$ such
that $K=\Q(\alpha)$ and $\Of=\Z[\alpha]$. In such a situation, there
is no obvious $\Phi(x)$ to work with.

Dedekind showed, however, that Under certain conditions we can still
make it work. Given a prime number $p\in\Z$, suppose we can find an
$\alpha$ such that $\Z[\alpha]\subset\Of$ has index not divisible by
$p$. Then factoring the minimal polynomial for $\alpha$ modulo $p$
gives the correct factorization of $(p)$ in $\Of$. This theorem was
announced by Dedekind in 1871 \cite{anzeige}; a proof appeared in 1878
\cite{Ded1878}; see \cite{GW2} for translations. It seems clear that
Kronecker was also aware of this fact, since he says he too started by
considering higher congruences.

This allowed one to hope, then, that an explicit factorization theory
could be based on a local approach: for each prime $p$, find a
generator $\alpha$ such that $p$ does not divide the index
$(\Of:\Z[\alpha])$. Then apply the theorem to find the
factorization. Dedekind says in \cite[\S4]{Ded1878} that he spent a
long time trying to prove that such an $\alpha$ always exists (see
\cite[p.~39]{GW2} for a translation).

Alas, this is not true: there exist number fields in which \emph{all}
of the indices have a common prime divisor. Dedekind pointed this out
(and stated the factorization theorem) in \cite{anzeige}, probably to
explain why he had needed to take a different route. Kronecker says in
his \emph{Grundz\"uge} \cite[\S25, p.~384]{Grundzuge} of 1882 that he
had found an example in 1858.

Both Dedekind and Kronecker pointed to this essential difficulty to
justify introducing a new approach: ideals in Dedekind's case, forms
in many variables in Kronecker's. Some years later, Zolotarev tried to
extend Kummer's theory directly in this style \cite{Zolo1}, but then
realized that his approach would fail for finitely many primes.
(Eventually, in a second paper \cite{Zolo2}, Zolotarev found still
another way to work around the difficulty.)  Dedekind's paper
\cite{Ded1878} was, as is clear from the introduction, prompted by
an announcement of Zolotarev's work. 

Kronecker also mentioned Zolotarev's attempt in
\cite[\S25]{Grundzuge}, where he stated the problem in terms of
discriminants. For each choice of $\alpha$, let
$d(\alpha)=\disc(\Phi(x))$ be the discriminant of its minimal
polynomial. Let $d_K$ be the field discriminant. Then
$d(\alpha)=m^2 d_K$, where $m$ is exactly the index
$(\Of:\Z[\alpha])$. Kronecker, who always preferred specific elements
to collections, thought about this as follows. The many element
discriminants $d(\alpha)$ have a common factor $d_K$ which is the
essential part, attached to the ``Gattungsbereich'' $K$ rather than to
a specific element. The other factors of $d(\alpha)$ (i.e., the
factors of $m$) are therefore ``inessential.'' So in the ``bad''
examples what is happening is that some prime $p$ is an inessential
divisor of every element discriminant. Such primes were the ``common
inessential discriminant divisors.''

The name is perhaps ill-chosen, because it is perfectly
possible\footnote{See footnote~\ref{ess_iness} on
  page~\pageref{ess_iness}.} for a prime $p$ to divide the
discriminant $d_K$ and \emph{also} divide the index
$(\Of:\Z[\alpha])$. Such a prime divisor is then both ``essential''
(it divides $d_K$) and ``inessential''! Dedekind's term ``index
divisor'' seems more appropriate. In the later literature, the ``index
$i(K)$ of the field $K$'' was defined to be the greatest common
divisor of the indices of all the generators of $K$; then Kronecker's
common inessential discriminant divisors are just the divisors of
$i(K)$. See \cite[2.2.1 item 3]{Nark} for information on more recent
work.

Kronecker's example ``in the thirteenth roots of 1'' is probably the
simplest one. He never gave the details, but they are probably as
Hensel gave them in his Ph.D.\ thesis \cite{HenselThesis} (see also
\cite[2.2]{Petri}). Let $\zeta$ be a primitive $13$-th root of
unity. There is a unique subfield $K$ of degree $4$ over
$\Q$.\footnote{This is global number field 4.0.2197.1 in
  \cite{LMFDB}.} Since the discriminant of $\Q(\zeta)$ is a power of
$13$, so is the discriminant of $K$ (in fact, $d_K=13^3$). It follows
from Kummer's work that the prime number $3$ is divisible by four
ideal primes in $K$, each of which has norm $3$; let $\fp$ be one of
these. Since $N(\fp)=3$, the field $\Of/\fp$ has three
elements. Consider some $\alpha\in K$. Since $K$ is a normal field,
the discriminant of the minimal polynomial of an integer in $K$ is the
square of the product of differences of the four roots, which are
integers in $K$.  Since there are only three congruence classes modulo
$\fp$, at least two of the roots must be congruent modulo $\fp$, i.e.,
one of these differences must be divisible by $\fp$. Since $\fp$ lies
above $3$, the discriminant $d(\alpha)\in\Z$ must be divisible by
$3$. Since $d_K$ is a power of $13$, the divisor $3$ is
inessential. This is true for any $\alpha$, so $3$ is a common
inessential discriminant divisor.

This set up the problem of determining exactly when this phenomenon
happens. One of the things that interests us about this problem is
that it was solved several times. Dedekind found a criterion in his
paper \cite{Ded1878}.\footnote{Hasse, in \cite[p.~456]{Hasse},
  attributes this criterion to Hensel and says it was the first
  success of Hensel's new methods, presumably meaning $p$-adic
  methods. In fact the criterion was first found by Dedekind and
  neither author used $p$-adic methods.} It is a sign of how little
Kronecker followed Dedekind's work that he suggested the problem of
common inessential discriminant divisors to Hensel for his Ph.D.\ in
1882. Hensel did not solve it completely in his
thesis,\footnote{Hensel later generalized the numerical condition in
  Kronecker's example to give a sufficient criterion for the existence
  of common inessential discriminant divisors, and even attempted to
  prove the condition was also necessary, which it is not. See the
  careful discussion in \cite[2.2]{Petri}.} but he published a
solution in 1894, in \cite{Hensel1894b}, which we translate
here. (Petri argues in \cite[2.4]{Petri} that the majority of the
results were known to Hensel before 1886.) In the first paragraph
gives the same criterion that had been found by Dedekind in
1878. While Hensel refers to Dedekind's paper, it is not clear how
carefully he had read it. In any case, he proceeds to find still
another criterion in the second half of the paper.

As Kronecker's student, Hensel does not work with ideals, but rather
with forms in several variables as in Kronecker's \emph{Grundz\"uge}
\cite{Grundzuge}. He probably learned this approach from Kronecker's
lectures, but those remained unpublished. Hensel refers to the
\emph{Grundz\"uge} as ``Kummer's \emph{Festschrift}'' because
Kronecker originally published it to commemorate the $50$th
anniversary of Kummer's doctorate.

We have tried, in our footnotes, to provide hints about how
Kronecker's approach works, without attempting a full account of
Kronecker's theory. For a modern attempt at explaining it, see
\cite{Weyl} or \cite{DivTheory}. The key thing to keep in mind for
this paper is that one considers a kind of ``generic algebraic
integer'': given an integral basis $\xi_1,\xi_2,\dots,\xi_n$
(Kronecker and Hensel call it a ``fundamental system'') one considers
the ``fundamental form''
\[ w_0=u_1\xi_1+ \dots +u_n\xi_n,\]
where the $u_1$ are indeterminates. This is a polynomial in the $n$
variables $u_1,u_2,\dots,u_n$; choosing integer values for the $u_i$
produces an algebraic integer in $K$.

Replacing the $\xi_i$ by their conjugates gives a conjugate $w_j$ of the
fundamental form; multiplying $x-w_j$ for all $j$ gives the
``fundamental equation'' for the domain. This is a polynomial in
$\Z[u_1,u_2,\dots,u_n,x]$. As before, specializing the $u_i$ to
integer values gives the polynomial of degree $n$ that has the
corresponding algebraic integer as a root.

Factoring the fundamental equation provides a method for finding the
factorization of a rational prime $p$ in the ring of integers of
$K$. This is just as in Dedekind's theorem: reduce the fundamental
equation modulo $p$ and factor it. The problem is that we are now
trying to factor a polynomial in $n+1$ variables.

The word ``Gattung'' means ``kind'' or ``genus''; Kronecker used it to
mean ``type of algebraic number.'' For example, $\sqrt2$ and
$3+\sqrt2$ are algebraic numbers of the same ``Gattung'' because each
is a rational function of the other. Algebraic numbers of the same
type belong to the same ``Gattungsbereich,'' which means something
like ``type domain.'' What we would call the base field Kronecker
called the ``domain of rationality.''  As usual, Kronecker did not
think in set-theoretic terms and would have avoided thinking of a
``Gattungsbereich'' as a completed whole.

Writing after Kronecker's death, Hensel seems a little more relaxed
about completed wholes---but also a little fuzzier. He seems to use
``Gattung'' and ``Gattungsbereich'' almost interchangeably for both a
field and its ring of integers. Since the word ``genus'' now means
something completely different we have opted to translate both words
as ``domain'' in most cases.

\section{The 1894 paper}

Hensel published two important papers in 1894. Both of them likely
contain material he originally submitted for his Habilitation in
1886. None of these were published at the time, and we know of them
only from Kronecker's notes. See \cite[2.3]{Petri} for a
reconstruction. The two 1894 papers were published after Kronecker's
death, perhaps because Hensel expected to find proofs for many of
these results among Kronecker's papers; see his explicit comment on
p.~\pageref{cannotfind} below.

The first paper published in 1894 was \cite{Hensel1894a}
``Untersuchung der Fundamentalgleichung einer Gattung f\"ur eine
reelle Primzahl als Modul und Bestimmung der Teiler ihrer
Discriminante'' (\textit{Journal f\"ur die Reine und Angewandte
  Mathematik} \textbf{113} (1894), 61--83). In it Hensel proved
something that had been stated by Kronecker in \cite{Grundzuge}: the
discriminant of the fundamental equation, which is a polynomial in
$\Z[u_1,u_2,\dots,u_n]$, has the discriminant of the field $K$ as its
largest integer factor. From this it follows that the factorization
modulo $p$ of the fundamental equation corresponds exactly to the
factorization of $p$.


The paper we translate builds on that to consider common inessential
discriminant divisors. Hensel wants to characterize when such divisors
occur. He finds several answers, the first of which is identical to
the one presented by Dedekind in 1878.

In our translation we have chosen to focus on getting the mathematical
content right, preserving Hensel's language, notations, and general
point of view. We have not tried (and would not have succeeded) to
preserve every nuance of meaning or to reproduce Hensel's grammar
precisely.

Our translation is based on the original publication in the
\textit{Journal f\"ur die Reine und Angewandte Mathematik}
(\textbf{113} (1894), 128--160); the original page numbers are
indicated in the margin. Hensel's own footnotes (which are few) are
marked by asterisks, while our annotations are given in numbered
footnotes. 

Hensel indicates theorems by using indented text; we have used
a modern simulacrum of the same device. Hensel does not signpost the
beginning or end of a proof; we have usually added such signposts in
the footnotes. We have retained Hensel's notation as much as
possible.

Every once in a while words have been inserted in square brackets when
we felt it would clarify the meaning. Hensel often uses $\gtrless$ to
indicate inequality or incongruence; we have silently substituted
$\neq$ or $\not\equiv$. We have also rendered ``ganze Functionen'' and
``ganze ganzzahlige Functionen'' as ``polynomials'' and ``integral
polynomials'' or ``polynomials with integer coefficients'',
respectively, without further comment. Hensel often says ``order''
when we would say ``degree''; he sometimes also uses ``dimension'' in
a similar sense. We have mostly translated ``degree'' when it was
unambiguous what was meant; see the footnotes.

\section*{Outline}

Hensel's paper contains five sections which he labels \S1 to \S5. The
main results in each section are as follows.

\S1. The main theorem here is that a prime $p$ is a common inessential
discriminant divisor in a field $K$ if and only if there are not
enough irreducible polynomials modulo $p$ to match the factorization
of $p$ in $K$. This result was also in \cite{Ded1878,GW2}.

\S2. Using the criterion Hensel just found seems to require knowing
the factorization of $p$, but in fact all we need to know is how many
prime divisors of $p$ in $K$ have a given degree. In this section
Hensel shows that one can determine this number without knowing the
factorization of $p$. Petri argues that the material in this paragraph
was not part of the Habilitation materials, hence was new in 1894.

\S3. The focus now changes to the index form
$\Delta(u_1,u_2,\dots,u_n)$ (with respect to a fixed integral basis).
Since $p$ is a common inessential discriminant divisor exactly when
plugging any $n$-tuple of integers into $\Delta$ results in a number
divisible by $p$. Hensel derives a general criterion to recognize when
a polynomial with integer coefficients has this property.

\S4. Kronecker had observed in \cite[\S25]{Grundzuge} an interesting
property of the index form in Dedekind's cubic field example. While
every value obtained by plugging integers into $\Delta(u_1,u_2,u_3)$
was divisible by $2$, there are integers from the cyclotomic field
$L=\Q(\zeta_3)$ for which we get values that are not divisible by
$3$. In this section Hensel shows that for any polynomial with integer
coefficients we can find an auxiliary field $L$ with this property.

\S5. Given the result in \S4, it is natural to ask which field we need
to use. Hensel shows that one can always choose a subfield of a
cyclotomic field of prime order.

\subsection*{Translation}

\vspace*{\baselineskip}

\begin{center}
  {\Large
  Arithmetical Investigations of the Common Inessential Discriminant
  Divisors of a Domain}\\[\baselineskip]
(by Mr.~K.~Hensel)\\
\rule{2in}{1pt}
\end{center}

\vspace{\baselineskip}

\begin{center}{\S1}\end{center}

\vspace{\baselineskip}

\marginpar{[128]} In a recently published work (this Journal,
vol. 111\footnote{This is a typo, as noted on page 160 of this issue;
  it should be volume 113. The paper is \cite{Hensel1894a}},
pp. 61--83) I considered the congruence of least degree modulo a prime
$p$ satisfied by the fundamental form\footnote{This alerts the reader
  that the $u_i$ in this equation are supposed to be
  indeterminates. Hensel, following Kronecker, uses ``form'' to mean a
  polynomial in several variables. The $\xi_i$ in this expression are
  what we would call an integral basis. One can think of $w_0$ as a
  ``generic algebraic integer.'' Multiplying $(w-w_0)$ with all its
  conjugates gives an equation of degree $n$ whose coefficients are in
  $\Z[u_1,u_2,\dots,u_n]$. This is the ``equation of smallest degree
  satisfied by $w_0$.'' Notice that $w_0$ is the fundamental form and
  $w$ is the variable in the polynomial it is a root of. Hensel will
  follow this notational pattern throughout.}
\[\tag{1.} w_0=u_1\xi_1+ \dots +u_n\xi_n\]
of a given domain\footnote{``Gattungsbereiches''} of the $n$-th
order.\footnote{We would say ``degree'' instead of ``order.''} The
main result, which will serve as the basis for this work, says that,
for any prime number, $w_0$ does not satisfy a congruence of degree
smaller than the degree $n$ of the domain.\footnote{In modern terms,
  the element $w_0$, considered modulo $p$, is integral of degree $n$
  over $\F_p[u_1,u_2,\dots,u_n]$.}

The congruence of lowest degree which is satisfied by $w_0$ modulo $p$
is made up in a simple manner from the congruences satisfied by $w_0$
modulo each of the prime divisors of $p$. Let $P$ be one of these
factors in the domain\footnote{``Bereich.'' We do not know why Hensel
  chooses $(\mathfrak{G})$ as the notation.} $(\mathfrak{G})$, and let
$\kappa$ be its degree.\footnote{``Ordnungzahl.'' We would call it the
  residual degree of $P$.}  Then $w_0$, with indeterminates
$u_1, \dots, u_n$, satisfies (modulo $P$) the congruence of degree
$\kappa$
\[\tag{2.} \mathfrak{F}(w)=w^\kappa + U^{(1)}(u_1 \dots
  u_n)w^{\kappa-1}+\cdots+U^{(\kappa)}(u_1 \dots u_n) \equiv 0
  \pmod{P},\] whose left side\footnote{Hensel thinks in terms of
  equations and congruences, not polynomials. The ``left hand side''
  of the congruence $\mathfrak{F}(w)\equiv 0$ is the minimal
  polynomial for $w_0~(\mathrm{mod} P$.}  is irreducible modulo the
prime $p$, while the coefficients are integral polynomials in
$u_1 \dots u_n$.

Let then
\[\tag{3.} p=P_1^{\delta_1}P_2^{\delta_2}\dots P_h^{\delta_h}\]
be the decomposition of $p$ into its prime factors in the domain
$(\mathfrak{G})$ and let
\[ \tag{4.}  \mathfrak{F}_1(w), \mathfrak{F}_2(w),\dots
  ,\mathfrak{F}_k(w) \] be the functions\footnote{In general, Hensel
  means ``polynomial'' when he says ``function.''} of lowest degree
having the fundamental form $w_0$ as a root modulo the $h$
\marginpar{[129]} corresponding distinct prime divisors
\[\tag{5.} P_1, P_2, \dots, P_h .\]
The polynomials have (as polynomials in $w$) degrees 
\[\kappa_1, \kappa_2,\dots, \kappa_h,\] 
which are equal to the degrees of the prime factors of
$p$.\footnote{So Hensel knows that the factorization of $p$ is
  determined by the factorization mod $p$ of the fundamental equation,
  which is the polynomial with coefficients in $\Z[u_1,u_2,\dots,u_n]$
  having $w_0$ as a root.} Then in the cited paper\footnote{This is
  \cite{Hensel1894a}.} (p.~75) it is shown that the congruence of
lowest degree satisfied by $w_0$ modulo the prime number $p$ will be
the following:
\[\tag{6.} \mathfrak{F}_1^{\delta_1}(w) \mathfrak{F}_2^{\delta_2}(w)\dots
\mathfrak{F}_h^{\delta_h}(w) \equiv 0 \pmod{p},\] 
and its degree in $w$ is
\[\tag{6$^{\mathrm a}$.}\kappa_1\delta_1+\kappa_2\delta_2+ \dots +
\kappa_h\delta_h=n.\] 
If we plug in $w_0$ for $w$, each of the $\mathfrak{F}_i(w_0)$ is
divisible by the divisor $P_i$, so the whole product is divisible by
$P_1^{\delta_1}\dotsc P_h^{\delta_h}$ and so divisible by $p$.

Instead of the fundamental form $w_0$ we now want to consider an
algebraic integer of the domain $(\mathfrak{G})$
\[\tag{7.}\xi_0 = a_1\xi_1+ a_2\xi_2+\dots + a_n\xi_n,\]
which is simply obtained\footnote{So $\xi_0$ is a particular algebraic
  integer, obtained as a linear combination of the integral basis
  $\xi_1,\xi_2,\dots, \xi_n$.} from $w_0$ by giving the unknowns
$(u_1, \dots, u_n)$ the integer values $(a_1, \dots, a_n)$. We would
like to investigate which congruence with integer coefficients is
satisfied by $\xi_0$ modulo $p$.

Obviously $\xi_0$ satisfies the congruence of $n$-th degree given by
(6.), since we are just replacing $(u_1, \dots, u_n)$ by
$(a_1, \dots, a_n)$. In general,\footnote{``In general'' here seems to
  mean that this notation will always be used.} if we have polynomials
\[\mathfrak{F}_1(w), \mathfrak{F}_2(w), \dots, \mathfrak{F}_h(w) \]
of degrees $\kappa_1, \kappa_2, \dots, \kappa_h$, we write the
corresponding\footnote{Hensel is introducing notation: these are the
  same functions as before, except that he has replaced the
  indeterminates $u_i$ by the integers $a_i$. The different variable
  $\xi$ indicates that this has been done. Notice that the functions
  $\mathfrak F_i(\xi)$ therefore depend on our choice of $\xi_0$,
  i.e., depend on the choice of the integers $a_i$. The main point of
  the notation is precisely not to have to show the dependence on the
  $u_i$ or the $a_i$.} functions
\[\mathfrak{F}_1(\xi), \mathfrak{F}_2(\xi), \dots,
\mathfrak{F}_h(\xi)\] (where the variable $w$ is now replaced by $\xi$
to indicate the difference). Then $\xi_0$ is a root of the congruence
of degree $n$ with integer\footnote{Whereas in (6.) the
  coefficients were polynomials in $n$ variables.} coefficients
\[\tag{8.} \mathfrak{F}_1^{\delta_1}(\xi) \dots,
\mathfrak{F}_h^{\delta_h}(\xi) \equiv 0 \pmod{p}.\] 
We know this to be the case because the individual numbers
$\mathfrak{F}_i(\xi_0)$ are obtained from the corresponding forms
$\mathfrak{F}_i(w_0)$ by replacing the unknowns $(u_1, \dots, u_n)$
\marginpar{[130]} by the integers $(a_1, \dots, a_n)$. Therefore each
of these are divisible a fortiori by the each of the prime factors as
before.

The fundamental form $w_0$, as proved in the aforementioned work, does
not satisfy any congruence whose degree is smaller than that of
(6.). But the congruence (8.) need not be the congruence of lowest
degree\footnote{For example, if we choose all of the $a_i$ divisible
  by $p$, the congruence of lowest degree will be the degree one
  polynomial $\xi$.} satisfied by $\xi_0$; for that to be
true we would need to choose $\xi_0$ appropriately.  We must then
investigate the following question:

\vspace{.5\baselineskip}

\hspace*{\fill}\parbox[t]{0.8\textwidth}{\hspace*{1em}Under
  what conditions will the algebraic integer $\xi_0$ satisfy no
  congruences modulo $p$ of degree less than $n$?}

\vspace{.5\baselineskip}

It is very easy to give a system of necessary conditions; we will
later prove that they are also sufficient. First,\footnote{Hensel
  doesn't do Lemmas, but this is one. Formally: if $\xi_0$ does not
  satisfy any congruence of degree less than $n$, then the polynomials
  appearing in (8.) must be irreducible modulo $p$.  Notice that these
  are the polynomials \emph{after} substituting the $u_i$ by the
  $a_i$, so they depend on the choice of $\xi_0$.} the $h$ integer
functions of $\xi$ in (8.), 
$\mathfrak{F}_1(\xi), \dots \mathfrak{F}_h(\xi)$ must be irreducible
modulo $p$. In fact,\footnote{Here begins the proof of the Lemma.} if
for example
\[\mathfrak{F}_1(\xi) \equiv F_1(\xi).G_1(\xi) \pmod{p},\]
where $F_1$ and $G_1$ are functions of $\xi$ of degree lower than
$\kappa_1$, then the same congruence also holds modulo the prime
divisor $P_1$ of $p$. Substituting $\xi$ by $\xi_0$ and observing
that $\mathfrak{F}_1(\xi_0)$ is divisible by $P_0$, we see the
congruence 
\[F_1(\xi_0).G_1(\xi_0) \equiv 0 \pmod{P_1},\]
from which it follows that one of these factors, say $F_1(\xi_0)$, is
divisible by the prime divisor $P_1$. Then $\xi_0$ clearly satisfies
the congruence modulo $p$ 
\[\tag{8$^{\mathrm a}$.} F_1^{\delta_1}(\xi)
  \mathfrak{F}_2^{\delta_2}(\xi) \dots \mathfrak{F}_h^{\delta_h}(\xi) 
  \equiv 0 \pmod{p},\] 
because the first factor on the left side is divisible by
$P_1^{\delta_1}$, while the others are divisible by
$P_2^{\delta_2}, \dots, P_h^{\delta_h}$. But the degree of
(8$^{\mathrm a}$.) is smaller than $n$, since it is
smaller\footnote{The contradiction ends the proof of the Lemma.} than
the degree of (8.).

Secondly, the $h$ irreducible polynomials $\mathfrak{F}_i(\xi)$ must
be distinct\footnote{Lemma 2, still under the running assumption that
  $\xi_0$ satisfies no congruence of degree lower than $n$.} modulo
$p$. If for example\footnote{Here begins the proof.}
$\mathfrak{F}_1(\xi)$ and $\mathfrak{F}_2(\xi)$ were congruent for
this modulus, then a fortiori
\[ \mathfrak{F}_1(\xi) \equiv \mathfrak{F}_2(\xi) \pmod{P_1P_2},\]
since $P_1P_2$ is a divisor of $p$. \marginpar{[131]}Substituting
again $\xi_0$ for $\xi$ and noticing that $\mathfrak{F}_1(\xi_0)$ is
divisible by $P_1$ and $\mathfrak{F}_2(\xi_0)$ is divisible by $P_2$,
it follows from the above congruence that $\mathfrak{F}_1(\xi_0)$ and
$\mathfrak{F}_2(\xi_0)$ are both divisible by the product
$(P_1P_2)$. Now if we take any two exponents $\delta_1,\delta_2$ with
$\delta_1\geq\delta_2$, the power $\mathfrak{F}_1^{\delta_1}(\xi)$
will, when $\xi = \xi_0$, be divisible by the product
$(P_1P_2)^{\delta_1}$, so a fortiori by
$P_1^{\delta_1}P_2^{\delta_2}$. Thus $\xi_0$ satisfies the congruence
\[\tag{8$^{\mathrm b}$.} \mathfrak{F}_1^{\delta_1}(\xi),
\mathfrak{F}_3^{\delta_3}(\xi) \dots \mathfrak{F}_h^{\delta_h}(\xi)
\equiv 0 \pmod{p},\] whose degree is smaller than that of (8.), so
smaller\footnote{End of the proof, again by contradiction.} than
$n$. So we have the following result:

\vspace{.5\baselineskip}

\noindent(A.)\hspace{\fill}\parbox{0.8\textwidth}{\hspace*{1em} If the
  number $\xi_0$ does not satisfy any congruence of degree less than
  $n$ modulo $p$, then the $h$ polynomials
  $\mathfrak{F}_1(\xi),\dots,\mathfrak{F}_h(\xi)$ in (8.) are all
  distinct and irreducible modulo $p$.}

\vspace{.5\baselineskip}

This theorem was given by Mr.~\emph{Dedekind} in his great
work\footnote{The original is ``grossen Arbeit.'' This is
  \cite{Ded1878}, which is really a short note rather than a
  full-length memoir, so ``grossen'' cannot mean ``large.'' For
  Dedekind it is an immediate consequence of his theorem about prime
  decomposition. See \cite{GW2}.} ``Ueber den Zusammenhang zwischen
der Theorie der Ideale und der Theorie der h\"oheren Congruenzen''
(\textit{Abh.\ der G\"ott.\ Gesellschaft} Volume 23), although in
slightly different form. It is noteworthy that he demonstrated that
for it to be possible to find such $h$ functions
$\mathfrak{F}_1(\xi), \dots, \mathfrak{F}_h(\xi)$, it is necessary and
sufficient that a number $\xi_0$ exists for which $p$ is not an
inessential divisor of the discriminant.\footnote{This is Theorem (IV)
  in \cite{Ded1878}; see \cite{GW2}. So at this point Hensel seems to
  know Dedekind's 1878 criterion for $p$ to be a common inessential
  discriminant divisor.} This result\footnote{``This result'' must be
  the necessary criterion in Theorem~A rather than the ``if and only
  if'' result just mentioned.}  enabled him to find a specific field
of degree three for which the prime number $2$ is a common inessential
divisor of all equation discriminants.  This alone shows that that the
theory of number fields\footnote{``Theorie der Gattungen.'' The same
  observation is made by Dedekind in the 1878 paper \cite{Ded1878} and
  by Kronecker in \cite[\S25]{Grundzuge}. This is exactly what
  Zolotarev tried to do but had to move beyond. Note the consensus
  here: the several different pioneers of algebraic number theory
  point to this issue to justify the complexity of their approaches.}
cannot be founded upon higher congruences. This\footnote{I.e., basing
  the whole theory on polynomial congruences.} can be done, however,
if, as in this and the previous work, we work with with the linear
\emph{form} $w_0=u_1\xi_1+\dots +u_n\xi_n$ with indeterminates
$u_1, \dots, u_n$, rather than a specific \emph{number} $\xi_0$ of the
domain. This is because in the previous work\footnote{The reference is
  to \cite{Hensel1894a}.} it was indeed established that in the
discriminant of the polynomial having $w_0$ as a root no prime $p$ is
contained as other than an essential divisor.\footnote{The theorem is
  that if you compute the discriminant of the fundamental equation you
  obtain a polynomial in $n$ variables $u_i$ whose content is $d_K$.}
The results on common inessential divisors of the discriminant that
follow in this paper have not, to my knowledge, been given
before.\footnote{In other words, Hensel is acknowledging that
  Theorem~A was proved by Dedekind but claims that his remaining
  theorems are new. This is not quite correct, since Dedekind also
  knew Theorem~B and knew that the criterion was sufficient, which
  Hensel will prove later; see Theorems~C and D.}

We now want to investigate when condition (A.) can actually
hold.\footnote{This section is about counting how many irreducible
  polynomials are available for Theorem~A. The first step is to group
  them by degree.} If the $h$ polynomials
$\mathfrak{F}_1(\xi)$, \dots, $\mathfrak{F}_h(\xi)$ are irreducible
\marginpar{[132]} modulo $p$, two of them, say $\mathfrak{F}_1(\xi)$
and $\mathfrak{F}_2(\xi)$, can only be congruent if they have the same
degree, i.e., if the degrees $\kappa_1$ and $\kappa_2$ of the
corresponding prime factors $P_1$ and $P_2$ are equal. So let us
arrange the distinct prime divisors $P_1,P_2,\dots,P_h$ by their
degrees $\kappa_1, \kappa_2,\dots \kappa_h$, grouping together those
that have equal degrees. Of the $h$ integers,
suppose that
\[ \text{there are }\lambda_1\text{ equal to }\kappa_1, \text{ then
}\lambda_2\text{ equal to }\kappa_2, \dots \text{ then
}\lambda_\gamma\text{ equal to }\kappa_\gamma,\] where $\lambda_1 +
\lambda_2 + \dots + \lambda_{\gamma}=h$ and the degrees $\kappa_1,
\kappa_2,\dots , \kappa_{\gamma}$ are all different. For the moment,
take $\kappa$ to be one of these $\gamma$ degrees\footnote{Now we
  focus on all the irreducible factors of a given degree.}
$\kappa_1, \dots ,\kappa_{\gamma}$ and let
\[\tag{9.} P^{(1)}, P^{(2)}, \dots, P^{(\lambda)}\]  
be the prime factors of $p$ whose degree is equal to $\kappa$.
Likewise, let
\[\tag{9$^{\mathrm a}$.} \mathfrak{F}^{(1)}(\xi), \mathfrak{F}^{(2 )},
  \dots, \mathfrak{F}^{(\lambda)}(\xi) \] be functions of degree
$\kappa$ satisfied by $\xi_0$ modulo the prime factors
$P^{(1)}, \dots P^{(\lambda)}$. We now want to know if it is possible
for these $\lambda$ functions to be irreducible and incongruent modulo
$p$.

Since $\mathfrak{F}^{(1)}(\xi)$ is irreducible [and has $\xi_0$ as a
root] modulo $P^{(1)}$, if $\xi_0$ satisfies another polynomial
congruence\footnote{The observation is that any polynomial such that
  $F(\xi_0)\equiv 0\pmod{P^{(i)}}$ must be divisible (modulo $p$) by
  the corresponding irreducible polynomial $\mathfrak{F}^{(i)}$.}
\[\varPhi(\xi) \equiv 0 \pmod{P^{(1)}},\]
then $\varPhi(\xi)$ must be divisible by $\mathfrak{F}^{(1)}(\xi)$
modulo $P^{(1)}$, because otherwise $\varPhi(\xi)$ and
$\mathfrak{F}^{(1)}(\xi)$ would have a greatest common divisor modulo
$P^{(1)}$, which contradicts by the irreducibility of
$\mathfrak{F}^{(1)}(\xi)$. Therefore for the modulus $P^{(1)}$ and
therefore\footnote{Since the polynomials all have integer
  coefficients.} for $p$ itself, we get a congruence of the form:
\[\varPhi(\xi) \equiv \mathfrak{F}^{(1)}(\xi)\varPhi^{(1)}(\xi) \pmod{p}.\] 
The same is true for the functions $\mathfrak{F}^{(2)}(\xi)$, \dots,
$\mathfrak{F}^{\lambda}(\xi)$ if they are also not decomposable for
$p$. 

Now\footnote{The residue field of each of the primes $P^{(i)}$ has
  $p^\kappa$ elements, all of which are roots of $x^{p^\kappa}-x$.}
all whole numbers from $(\mathfrak{G})$, and hence also $\xi_0$,
satisfy the congruence
\[\xi^{p^{\kappa}}-\xi \equiv 0 \pmod{P^{(i)}} \hspace{1in}
  \textit{\scriptsize $(i=1,2,\dots,
    \lambda)$},\] for each of the $\lambda$ divisors of degree
$\kappa$ $P^{(1)}, \dots, P^{(\lambda)}$.  Hence the expression
$(\xi^{p^{\kappa}}-\xi)$ is divisible modulo $p$ by each of the
$\lambda$ functions $\mathfrak{F}^{(1)}(\xi),$ \dots
$\mathfrak{F}^{\lambda}(\xi)$, if they are assumed to be
irreducible. \marginpar{[133]} If those functions are incongruent
modulo $p$, then the function $(\xi^{p^{\kappa}}-\xi)$ must be
divisible by their product, and so it must contain at least $\lambda$
irreducible factors modulo $p$ of degree $\kappa$.\footnote{Hensel has
  shown, then, that any irreducible polynomial of degree $\kappa$ is a
  divisor of $\xi^{p^\kappa}-\xi$ modulo $p$. In modern terms,
adjoining a root of an irreducible polynomial of degree $\kappa$ to
$\F_p$ always gives the same field, namely the splitting field of
$\xi^{p^\kappa}-\xi$. In fact, Hensel also needs to know that any
polynomial whose degree \emph{divides} $\kappa$ is a factor modulo $p$
of $\xi^{p^k}-\xi$. For that he quotes his older paper
\cite{HenselFF}, which is one of many nineteenth century papers
dealing with ``higher congruences'' that we would describe as being
about the theory of finite fields.}  If we consider all the
irreducible factors modulo $p$ one finds\footnoteB{Compare with my
  paper: Untersuchung der ganzen algebraischen Zahlen eines
  Gattungsbereiches f\"ur einen beliebigen algebraischen Primdivisor;
  this Journal, volume 101, pages 140 and 141.}  that
$(\xi^{p^{\kappa}}-\xi)$ is the product of all irreducible polynomials
whose degree is equal either to $\kappa$ or to a divisor of $\kappa$
and that $(\xi^{p^{\kappa}}-xi)$ has exactly
\[\tag{10.} \bar{g}(\kappa)=\frac{1}{\kappa}(p^{\kappa} - \sum
  p^{\frac{\kappa}{q}} + \sum p^{\frac{\kappa}{qq'}} - \sum
  p^{\frac{\kappa}{qq'q''}} + \dots)\] distinct irreducible divisors
of degree $\kappa$, where $q, q', q'', \dots$ are the distinct prime
factors of $\kappa$.\footnote{The point is that we have a polynomial
  of degree $p^\kappa$ which is the product of all irreducible
  polynomials mod $p$ whose degree divides $\kappa$. Writing
  $p^\kappa$ as the sum of those degrees and using M\"obius inversion
  gives formula (10.) for the total number of distinct irreducible
  polynomials of degree $\kappa$ in $\F_p[\xi]$. The same formula is also
  found in Dedekind's \emph{Abri\ss} \cite{Abriss}, but Dedekind does
  not quote it in his 1878 paper.}  So if
$\lambda > \bar{g}(\kappa)$ then it is not possible for the $\lambda$
irreducible functions $\mathfrak{F}^{(1)}(\xi)$, \dots
$\mathfrak{F}^{\lambda}(\xi)$ to be distinct modulo $p$. If we now
apply this result to all $\gamma$ of the distinct degrees
$\kappa_1, \dots ,\kappa_{\lambda}$ of the prime factors of $p$, we
obtain from theorem (A.):\footnote{In the statement of this theorem
  Hensel uses ``real prime'' to refer to a prime in $\Z$. Similarly, he
  later uses ``real integer'' for an element of $\Z$. The modern usage
  is ``rational integer'' and ``rational prime,'' but we have
  preserved Hensel's words.}

\vspace{.5\baselineskip}

\noindent (B.\footnote{This is just theorem~A plus an explicit count
  of the number of irreducible polynomials in degree $\kappa$ in
  $\F_p[\xi]$.})\hspace{\fill}\parbox{.8\textwidth}{\hspace*{1em}Suppose that
  $p=P_1^{\delta_1} \dots P_h^{\delta_h}$ is the decomposition of a
  real prime number $p$ in a domain
  $(\mathfrak{G})$, and that among the $h$ nonequivalent prime factors
  $P_1, P_2, \dots, P_h$ there are
  \[ \lambda_1\text{ of degree }\kappa_1,\]
  \[ \lambda_2\text{ of degree }\kappa_2,\]
  \[\vdots\]
  \[ \lambda_\gamma\text{ of degree }\kappa_\gamma.\] So we can find a
  number $\xi_0$ in the domain $(\mathfrak{G})$ that does not satisfy
  any congruence modulo $p$ of degree less than $n$ only if
  \[\tag{11.}\lambda_1 \leq \bar{g}(\kappa_1),\quad \lambda_2 \leq
  \bar{g}(\kappa_2),\quad \dots, \quad \lambda_{\gamma} \leq
  \bar{g}(\kappa_{\gamma})\] 
  holds (where the $\gamma$ whole numbers $\bar{g}(\kappa)$ are as in
  (10.)). If however, even one of these conditions is not met, then
  \emph{every} number $\xi_0$ from $(\mathfrak{G})$ satisfies a 
  polynomial congruence modulo $p$ of degree less than $n$.}

\vspace{.5\baselineskip}

\marginpar{[134]} It should now be proved that condition (11.) is also
sufficient\footnote{So another proof is beginning here: that the
  conditions (11.) imply the existence of at least one $\xi_0$ with
  the desired property. This is Theorem (IV) in
  \cite{Ded1878}. Hensel's proof is identical to Dedekind's. See
  \cite{GW2}.} to guarantee that at least one number $\xi_0$ from
$(\mathfrak{G})$ satisfies no congruence modulo $p$ of degree less
than $n$.  Let $P$ be one of the $h$ prime divisors of $p$ and let
$\kappa$ be its degree. Then we can choose\footnote{This is the lemma
  to be proved next. Given an irreducible polynomial of degree
  $\kappa$ in $\F_p[x]$, we can choose $\xi_0$ so that it is a root of
  that polynomial modulo $P$.} $\xi_0$ to satisfy the irreducible
polynomial congruence
\[\mathfrak{F}(\xi)\equiv 0 \pmod{P},\]
where $\mathfrak{F}(\xi)$ is one of the $\bar{g}(\kappa)$ irreducible
divisors of degree $\kappa$ of $(\xi^{p^{\kappa}}-\xi)$.

In\footnote{Here starts the proof of the lemma.} the congruence
\[\xi^{p^{\kappa}}-\xi \equiv \mathfrak{F}(\xi)\varPhi(\xi)\pmod{P},\]
both the left and ride side disappear for as many incongruent values
of $\xi$ as the degree (namely for the $p^\kappa$ congruence classes
modulo $P$ of numbers in the domain $(\mathfrak{G})$). There must
therefore exist a number $\xi_0$ for which $\mathfrak{F}(\xi_0)$ is
divisible by $P$, because if not the function $\varPhi(\xi)$ of degree
$(p^\kappa-\kappa)$ would vanish modulo $P$ for $p^\kappa$ incongruent
values of $\xi$, which is not possible. We can
therefore\footnote{Lemma has been proved.} choose $\xi_0$ to be a root
of the chosen irreducible congruence of degree $\kappa$
\[\mathfrak{F}(\xi) \equiv 0 \pmod{P}.\]

Now choose\footnote{We have shown that for each prime $P$ of degree
  $\kappa$ and each irreducible polynomial of degree $\kappa$ we can
  find an integer that is a root of that polynomial modulo $P$. Now we
  apply this to each of the prime factors of $p$.} $h$ irreducible
functions $\mathfrak{F}_1(\xi)$, $\mathfrak{F}_2(\xi)$, \dots,
$\mathfrak{F}_h(\xi)$, all incongruent modulo $p$, whose degrees are
respectively equal to $\kappa_1,\kappa_2,\dots\kappa_h$, so that each
of the $\mathfrak{F}_i(\xi)$ is a divisor of $\xi^{p^{x_i}}-\xi$. Such
a system of $h$ functions can only exist when condition (11.) is
satisfied. When it is, we can find such a system, since the quantity
of irreducible functions incongruent modulo $p$ of degrees
$\kappa_1,\dots,\kappa_h$ is larger than the number of functions we
need. Let then
\[\tag{12.} \xi_0^{(1)}, \xi_0^{(2)}, \dots, \xi_0^{(h)}\]
be $h$ integers from $(\mathfrak{G})$ chosen so that for each $i$ the
integer $\xi_0^{(i)}$ satisfies modulo $P_i$ the congruence 
\[\tag{12$^{\mathrm a}$.}\mathfrak{F}_i(\xi_0^{(i)}) \equiv 0
\pmod{P_i} \hspace{1in} (i=1, \dots, h).\] As we proved above, we can
find $h$ such numbers $\xi_0^{(i)}$.

Moreover, we can also
assume\footnote{Another little lemma.} from the outset that the left
side in (12$^{\mathrm a}$.) is not divisible by $P_i^2$, \marginpar{[135]}
so that
\[\tag{12$^{\mathrm b}$.}\mathfrak{F}_i(\xi_0^{(i)}) \not \equiv 0
\pmod{P_i^2} \hspace{1in} (i=1,2, \dots h).\] 
This is possible because if for some $i$ 
\[\mathfrak{F}_1(\xi_0^{(1)}) \equiv 0 \pmod{P_1^2},\]
one can substitute $\xi_0^{(1)}$ by
$\bar\xi_0^{(1)}=(\xi_0^{(1)}+\pi_1)$, where $\pi_1$ is a whole number
divisible by $P_1$ but not $P_1^2$.  Then using\footnote{Hensel
  doesn't discuss the denominators in the Taylor expansion. What is
  actually needed is a version of Taylor's theorem for polynomials:
  express the formal polynomial $f(X+Y)$ as a polynomial in $Y$ with
  coefficients in $\Z[X]$.} \emph{Taylor}'s theorem, we have:
\[\mathfrak{F}_1(\bar\xi_0^{(1)})= \mathfrak{F}_1(\xi_0^{(1)}+\pi_1)=
  \mathfrak{F}_1(\xi_0^{(1)})+\pi_1\mathfrak{F}_1^{\,\prime}(\xi_0^{(1)})
  +\frac{1}{2}\pi_i^2\mathfrak{F}_1^{\,\prime\prime}(\xi_0^{(1)})+
  \dots.\] According to our assumptions, the first term as well as the
third and all subsequent terms\footnote{If $p=2$ it does not seem
  clear that the third term is divisible by $P_1^2$. See the previous
  footnote.} are divisible by $P_1^2$. But in the second term
$\mathfrak{F}_1^{\,\prime}(\xi_0^{(1)})$ is not divisible by $P_1$
(because $\mathfrak{F}_1^{\,\prime}(\xi)$ cannot have a common
divisor\footnote{Hensel is assuming, perhaps without noticing it, that
  he doesn't have to worry about the possibility that
  $\mathfrak F_i^{\,\prime}(\xi)\equiv0$. He is correct because finite
  fields are perfect.} with the irreducible polynomial
$\mathfrak{F}_1(\xi)),$ we see the following congruence:
\[\mathfrak{F}_1(\bar \xi_0^{(1)}) \equiv
  \pi_1\mathfrak{F}_1^{\,\prime}(\xi_0^{(1)}) \not \equiv 0
  \pmod{P_1^2}.\] Thus,\footnote{The lemma has been proved.} the $h$
numbers $\xi_0^{(1)}, \dots, \xi_0^{(\lambda)}$ can be chosen from the
beginning in such a way that the $h$ terms
\[\mathfrak{F}_1(\xi_0^{(1)}), \mathfrak{F}_2(\xi_0^{(2)}), \dots ,
\mathfrak{F}_h(\xi_0^{(h)})\] 
are divisible once and only once by the corresponding divisors
\[P_1, P_2, \dots, P_h.\]

If this occurs, then it is also possible to find\footnote{In the
  parallel passage of \cite{Ded1878}, Dedekind simply invokes the
  Chinese Remainder Theorem from \cite{SuppX}; Hensel is going to give
  a proof.} an algebraic number $\xi_0$, so that
\[\tag{13.} \xi_0 \equiv \xi_0^{(i)} \pmod{P_i^2} \hspace{1in} (i=1,2,
  \dots, h).\]
There always exists\footnote{Here begins the proof of the Chinese
  Remainder Theorem in this situation.} a number $\varPi_1$ in the
domain $(\mathfrak{G})$ that is not divisible by $P_1$, but is
divisible by all other divisors of $p$. We can find another number
$\varrho_1$ so that\footnoteB{Set $\varrho_1=x_1+\pi_1y_1$, where
  $\pi_1$ is divisible by $P_1$ exactly once. Then the algebraic
  numbers $x_1$ and $y_1$ are defined by the system of linear
  congruences for the Modulus $P_1$
    \[ x_1\varPi_1^2 \equiv 1 \pmod{P_1},\qquad
    \frac{(x_1+\pi_1y_1)\varPi_1^2-1}{\pi_1} \equiv 0 \pmod{P_1},\] which
    always has a solution.}  
\[ \varepsilon_1  = \varrho_1\varPi_1^2 \equiv 1 \pmod{P_1^2}. \]
Then the number $\varepsilon_1$ is divisible by each of the divisors
$P_2^2$, \dots $P_h^2$, \marginpar{[136]}
while the remainder of division by $P_1^2$ is
$1$. So we take 
\[\varepsilon_1, \varepsilon_2, \dots, \varepsilon_h\]
to be $h$  algebraic integers chosen so that
\begin{align*}
  \varepsilon_i &\equiv 1 \pmod{P_i^2}\\
  \varepsilon_i &\equiv 0 \pmod{P_l^2}, \hspace{1in} (i \neq l)
\end{align*}
and set
\[\tag{14.} \xi_0=\epsilon_1\xi^{(1)}+\epsilon_2\xi_0^{(2)}+
\dots+\epsilon_h\xi_0^{(h)}.\] 
Then $\xi_0$ satisfies the $h$ conditions (13.).\footnote{So we have
  proved the Chinese Remainder Theorem.} Now from (12$^{\mathrm a}$) and
(12$^{\mathrm b}$), we have  
\[\tag{15.}
\begin{cases}
  \mathfrak{F}_i(\xi_0) \equiv \mathfrak{F}_i(\xi_0^{(i)}) \equiv 0
  \pmod{P_i}\\ \rule{0ex}{3ex} \mathfrak{F}_i(\xi_0) \equiv
  \mathfrak{F}_i(\xi_0^{(i)}) \not \equiv 0 \pmod{P_i^2}\end{cases}.\]
Moreover, the expression $\mathfrak{F}_i(\xi_0)$ is not divisible by
any prime divisors $P_l$ distinct from $P_i$. (If so, since
$\mathfrak{F}_l(\xi_0)$ is certainly divisible by $P_l$, the
irreducible polynomials $\mathfrak{F}_i(\xi)$ and
$\mathfrak{F}_l(\xi)$ would have a common divisor, which
means\footnote{Because they are irreducible.} they would be congruent
to each other, which conflicts with the above assumption.)

Since $\xi_0$ is chosen according to the condition (14.), for each $i$
this number satisfies the irreducible congruence of degree $\kappa_i$ 
\[\mathfrak{F}_i(\xi_0) \equiv 0 \pmod{P_i}\] 
and this expression is not divisible by any other\footnote{Hensel
  means ``by any other prime factor of $p$.''} prime
factor, which means:\footnote{Here $\sim$ denotes equivalence of
  divisors in Kronecker's sense. We would say that $P_i$ is the
  greatest common divisor of $p$ and $\mathfrak F_i(\xi_0)$.} 
\[p+u. \mathfrak{F}_i(\xi_0) \sim P_i.\]
If this is the case, then we can prove exactly the same way as was
done for the fundamental form of $w_0$ in \S3 of the previous
work\footnote{Again, this is \cite{Hensel1894a}.} that the congruence of
degree $n$ 
\[\mathfrak{F}_1^{\delta_1}(\xi)\mathfrak{F}_2^{\delta_2}(\xi) \dots
  \mathfrak{F}_h^{\delta_h}(\xi) \equiv 0 \pmod{p},\] is the
smallest\footnote{I.e., is of the smallest degree.} that $\xi_0$
satisfies, which means it is the element we need for the
converse\footnote{So this concludes the proof of the converse: if the
  inequalities (11.) hold, then an element $\xi_0$ as in the theorem
  can be found. See \cite{GW2} for a numerical example.} of
Theorem~(B.).

Recalling that the degree $\kappa$ of a prime divisor $P$ of $p$
coincides with the degree (as a polynomial in $w$) of the
corresponding irreducible modulo $p$ factor $\mathfrak F(w)$ from
(2.), we can restate\footnote{So here we are back to the factorization
  of the fundamental equation, which is a polynomial in $w$ with
  coefficients in $\Z[u_1,u_2,\dots,u_n]$.} the result without
using\footnote{The idea, of course, is that in general it is hard to
  find the factorization of $p$, and especially so when the condition
  in the theorem holds. Alas, factoring the fundamental equation is
  hard as well. Hensel will address this in the next section.}
the prime factorization of $p$, in the following manner:

\vspace{.5\baselineskip}
\marginpar{[137]}

\noindent(C.)\hspace{\fill}\parbox{0.8\textwidth}{If
\[\mathfrak{F}(w) \equiv \mathfrak{F_1}^{\delta_1}(w) \dots
\mathfrak{F}_h^{\delta_h}(w) \pmod{p}\] 
is the decomposition of the fundamental equation of a field
$(\mathfrak{G})$ into its irreducible factors modulo $p$, and if the 
numbers  
\[\lambda_1, \lambda_2, \dots, \lambda_{\gamma},\]
indicate how many of the $h$ factors $\mathfrak{F}_1(w), \dots,
\mathfrak{F}_h(w)$ have corresponding degree 
\[\kappa_1, \kappa_2, \dots, \kappa_{\gamma},\]
then \emph{all} values $\xi_0$ of domain modulo satisfy a congruence
of degree lower than $n$ if and only one of the $\gamma$ conditions
\[\lambda_i > \bar{g}(\kappa_i) (i=1,2, \dots, \gamma)\]
is satisfied. Here the term $\bar{g}(\kappa)$ is
\[\bar{g}(\kappa)=\frac{1}{\kappa}(p^\kappa - \sum
p^{\frac{\kappa}{q}}+\sum p^{\frac{\kappa}{qq'}}-\dots),\] 
and $q, q', \dots $ are the distinct prime factors of the number
$\kappa$.} 

\vspace{.5\baselineskip}

Now\footnote{The next few paragraphs relate the field discriminant and
  the discriminant of an element, introducing the notions of ``index''
  and of ``inessential divisor.'' Hensel is still tracking Dedekind
  quite closely, but he would have seen some of this in Kronecker as
  well. The first step is to relate the fact that $\xi_0$ satisfies a
  congruence of degree less than $n$ to its index.}  take $\xi_0$ to
be some number in the domain $(\mathfrak{G})$; listing the first $n$
powers of $\xi_0$ in terms of the fundamental
system\footnote{``Fundamental system'' is Kronecker's name for an
  integral basis.} $\xi_1,$ \dots, $\xi_n$, we get $n$ equations with
rational integer coefficients
\[\tag{16.}\left\{
\begin{aligned}
  1 &= a_{10}\xi_1 + \dots +a_{n0}\xi_n,\\
  \xi_0 &= a_{11}\xi_1 + \dots +a_{n1}\xi_n,\\
  \vdots & \\
  \xi_0^{n-1}&= a_{1,n-1}\xi_1 + \dots +a_{n,n-1}\xi_n.
\end{aligned}
\right. \]
Now $\xi_0$ satisfies a polynomial congruence modulo $p$ of
degree less than $n$ if and only if the determinant 
\[|a_{ik}|\qquad\qquad
\genfrac{}{}{0pt}{}{(i=1, \dots, n)}{\qquad\quad\!(k=i, 1, \dots,
  n-1)}\]
of the $n$ linear equations (16.)  is divisible by $p$. This is
because\footnote{A bit of linear algebra modulo $p$: we want the
  system $A[\xi]= [0]$ to have a nontrivial solution mod $p$, which
  requires the determinant to be zero mod $p$.} only in this case can
we find $n$ numbers $A_0,$ $A_1,$ \dots , $A_{n-1}$ (not all divisible
by $p$), such that the sum of the first equation multiplied with
$A_0$, the second with $A_1,$ \dots, and the last with $A_{n-1}$,
gives
\[A_0+A_1\xi_0+ \dots + A_{n-1}\xi_0^{n-1} \equiv 0 \pmod{p}\]
by making the \marginpar{[138]} coefficients of $\xi_1, \dots , \xi_n$
on the right side of the equation all\footnote{Since
  $\xi_1$,\dots,$\xi_n$ is an integral basis, the only way an
  algebraic integer will be divisible by $p$ is by having all the
  coefficients divisible by $p$.}  divisible by $p$. When we form the
$n$ systems of equations from (16.)  and consider the $n$
conjugate\footnote{This is basically matrix multiplication: the matrix
  whose columns are the powers of $\xi_0$ and its conjugates is equal
  to $[a_{ik}]$ times the matrix whose columns are the integral basis
  and its conjugates. Hensel thinks of $n$ conjugate domains rather
  than doing the computation in a normal closure.} domains to
${\mathfrak{G}}$, then we see the validity of the equation
\[\mathfrak{D}(\xi_0)=|a_{ik}|^2.D,\]
where $\mathfrak{D}(\xi_0)$ is the discriminant of the equation for
$\xi_0$ and $D$ is the discriminant of the domain $(\mathfrak{G})$. 

Each equation discriminant consists, then, of two essentially
different parts: on the one hand, the domain discriminant $D$, which
is the same in all discriminants, and on the other the squared
determinant $|a_{ik}|^2$, which is dependent on the choice of
$\xi_0$. For this reason, \emph{Kronecker} called the first the
essential and the second the inessential divisor of the discriminant
$\mathfrak{D}(\xi_0)$. A prime $p$ is contained in $|a_{ik}|^2$ (and
so is an inessential divisor of the discriminant
$\mathfrak{D}(\xi_0)$) if and only if $\xi_0$ satisfies a congruence
modulo $p$ of degree less than $n$.\footnote{Dedekind \cite{Ded1878}
  called the (absolute value of the) determinant $|a_{ik}|$ the
  \emph{index} of the algebraic integer $\xi_0$. Hensel seems to be
  content not to have a name for it.}  From our Theorem~(C.), it
follows now that the first part $|a_{ik}|^2$ of the discriminant
(although it depends on $\xi_0$) can contain factors which remain the
same whatever $\xi_0$ is chosen, and so cannot be removed by an
appropriate choice of $\xi_0$. These ``common inessential divisors''
of all equation discriminants of a domain are the primes $p$ (and only
these) for which every number $\xi_0$ of the domain satisfies a
congruence of lower than $n$-th degree. By applying Theorem~(C.) we
get the following:

\vspace{\baselineskip}

\noindent(D.)\hspace{\fill}\parbox{0.8\textwidth}{
\hspace*{1em}If
\[F(w) \equiv \mathfrak{F}_1^{\delta_1}(w) \dots
  \mathfrak{F}_h^{\delta_h}(w)\pmod{p}\] is the decomposition of the
fundamental equation of a field $(\mathfrak{G})$ into its irreducible
factors modulo $p$, and if the numbers
$\lambda_1, \dots , \lambda_{\gamma}$ indicate how many factors have
degree $\kappa_1, \dots , \kappa_{\gamma}$ as polynomials in $w$, then
$p$ is a common inessential divisor of all equation discriminants
$\mathfrak{D}(\xi_0)$ from $(\mathfrak{G})$ if and only if at least
one of the $\gamma$ conditions
\[\lambda_1> \bar{g}(\kappa_1), \dots,
\lambda_{\gamma}>\bar{g}(\kappa_\gamma)\] 
is satisfied.}

\vspace{\baselineskip}

\begin{center} \S2\marginpar{[139]} \end{center}

\vspace{\baselineskip}

The result from the previous section can also be expressed in another
form that is remarkable in that to apply it we do not need to know the
decomposition of $p$ within the domain $(\mathfrak{G})$ or the
factorization of the fundamental equation modulo $p$.\footnote{The
  fundamental objection to the criterion above is that in order to use
  it we need to know how $p$ factors, or, equivalently, we need to be
  able to factor the fundamental equation modulo $p$. This can be hard
  to do in general. What Hensel notes in this section is that in fact
  one does not need to know the full factorization. It suffices to
  know, for each $\kappa=1,2,\dots,n$, the number $\lambda_\kappa$ of
  distinct primes of degree $\kappa$ occurring in the factorization of
  $p$. Then we can compare that with the number $\bar{g}(\kappa)$ of
  irreducible polynomials of degree $\kappa$ in $\F_p[x]$ to determine
  whether $p$ is a common inessential discriminant divisor. The goal
  of this section is to present a way of computing $\lambda_\kappa$
  without finding the full factorization of $p$. From this point on
  Hensel is venturing beyond Dedekind's 1878 paper \cite{Ded1878,GW2}.}

Let $P$ be an arbitrary prime factor of $p$ and let $\kappa$ be its
[residual] degree. Now if 
\[w_0=u_1\xi_1+ \dots+u_n\xi_n\] is a fundamental form\footnote{That
  is, $\{\xi_1,\xi_2,\dots,\xi_n\}$ is an integral basis and the $u_i$
  are indeterminates.} for the domain $(\mathfrak{G})$, we set
\[w_h=u_1\xi_1^{p^h}+ \dots+u_n\xi_n^{p^h}\qquad\qquad (h=0,1,
  \dots).\] We know from the previous work\footnote{Namely
  \cite{Hensel1894a}.} (page~65), that only the first $\kappa$ of these
infinitely many forms\footnote{Most readers will not be used to
  Kronecker's form-based version of algebraic number theory. The thing
  to note is that the form $w_h$ is a ``generic element'' of the ideal
  generated by $\xi_i^{p^h}$. Whenever Hensel talks of such a form,
  one can translate to Dedekindian terms by considering the ideal
  generated by the coefficients.

  The forms $w_h$ are lifts of the images of $w_0$ under the Frobenius
  automorphism modulo $P$, which is of order $\kappa$. This gives the
  congruence claims that follow immediately. Hensel does not
  have any of this language at his disposal, of course.}
\[w_0, w_1, \dots, w_{\kappa-1}\]
are distinct modulo $P$ (for indeterminate $u_1, \dots , u_n$). In
fact, $w_\kappa \equiv w_0$, $w_{\kappa+1} \equiv w_1, \dots$, and in
general 
\[w_{h+\kappa} \equiv w_h \pmod{P}.\] From this it follows that the
linear form\footnote{This is the key object for this section. In
  Dedekindian terms, we are considering the ideal $I_\nu$ generated by
  the elements $\xi_i^{p^\nu}-\xi_i$. Notice that these ideals depend
  on the choice of integral basis.}
\[w_{\nu}-w_0= u_1(\xi_1^{p^{\nu}}-\xi_1)+u_2(\xi_2^{p^{\nu}}-\xi_2)+
  \dots +u_n(\xi_n^{p^{\nu}}-\xi_n)\] is divisible by the prime
divisor $P$ if and only if the number $\nu$ is a multiple of the
degree $\kappa$. If this is the case, then the following simple
considerations show that $P$ is only contained once in that linear
form.\footnote{The claim is that when $\kappa|\nu$ the ideal $I_\nu$
  is divisible exactly once by each prime of degree $\kappa$ appearing
  in the factorization of $p$. It is clear that $P$ divides
  $w_\nu-w_0$, so the key thing to prove is that $P^2$ does not. The
  next two paragraphs appear to provide a proof of this claim. In
  fact, however, what they show is that the forms $w_{\nu}-w_0$ may
  need to be modified so that this is true. Hensel phrases this as
  modifying the integral basis, but after the modification he suggests
  the list $\xi_1,\xi_2,\dots,\xi_n$ is no longer an integral basis.

  The ``proof'' is divided into two cases: when $P^2$ divides $p$ and
  when it does not. In the first case, we have an actual proof. In the
  second a change to the integral basis is needed. Hensel is assuming
  we do not know the factorization of $p$, however, so he will make
  the modification in any case.}

If, first, $P$ is a multiple divisor of $p$, then $P^2$ is contained
in $p$ and we know from \textit{Fermat}'s [little] theorem, that for
all integer values of $u_1, \dots , u_n$ the congruence
\[w_{\nu} \equiv w_0^{p^\nu} \pmod{p}\]
holds. The same is also fulfilled a fortiori modulo the divisor $P^2$
of $p$. If it were true that even for indeterminate $u_1, \dots , u_n$
we had $P^2$ dividing the linear form $w_{\nu}-w_0$, then
\emph{all} numbers $\xi_0=a_1\xi_1+ \dots + a_n\xi_n$ in
$(\mathfrak{G})$ would satisfy the congruence
\[\xi_0^{p^{\nu}}-\xi_0 \equiv 0 \pmod{P^2}.\]
That this is not the case\footnote{So when $P^2$ divides $p$ the form
  $w_{\nu}-w$ is never divisible by $P^2$. The argument involves
  choosing a uniformizer at $p$.} can be easily seen if we
consider, for example, \marginpar{[140]} $\xi_0=\pi$, where $\pi$ is
divisible by $P$, but not by $P^2$; in this case, the left side of the
congruence reduces modulo $P^2$ to $(-\pi)$, since $\pi^{p^\nu}$ is
clearly divisible by $P^2$.


Next, if $P$ divides $p$ only once and if the linear
form $w_{\nu}-w_0$ for indeterminate $(u_1, \dots, u_n)$ is divisible
by $P^2$, then we must have 
\[\xi_i^{p^{\nu}}-\xi_i \equiv 0 \pmod{P^2}\qquad\qquad (i=1,2, \dots,n)\]
for the $n$ elements of the fundamental system, which as can easily be
seen is generally\footnote{In the unramified case it is indeed possible for
  $\xi_i^{p^\nu}-\xi_i$ to be divisible by $P^2$. A simple example is
  to take $K=\Q(\sqrt3)$ with $\xi_1=1$ and $\xi_2=\sqrt3$. Let
  $p=11$, Then $\xi_1^{11}-\xi_1=0$ and
  \[\xi_2^{11}-\xi_2 = (\sqrt3)^{11}-\sqrt3 = 242\sqrt3 =
  2\cdot11^2\sqrt3.\] If, as Hensel suggests, we instead use
  $\xi_1=12$, then it works, since $12^{11}-12$ is divisible by $11$
  only once.

  The most dramatic example is the cyclotomic field generated by an
  $\ell$-th root of unity when $p\equiv 1\pmod{\ell}$. If we take the
  standard integral basis, then $w_\nu=w_0$ for all $\nu$, and
  $w_\nu-w_0=0$.

  For a cubic example, let $K=\Q(\alpha)$ with
  $\alpha^3 - 6\alpha^2 - 9\alpha - 1=0$ (number field 3.3.3969.2 in
  \cite{LMFDB}). The integral basis is $(1,\alpha,\alpha^2)$ and when
  $p=5$ both $\alpha^5-\alpha$ and $\alpha^{10}-\alpha^2$ turn out to
  be divisible by the square of one of the primes dividing $5$.} not
the case. For if the system $(\xi_1, \dots, \xi_n)$ did have this
property, we could modify it, without changing its
character\footnote{Hensel is correct that this does not change the
  reduction mod $p$ of the $\xi_i$, but of course they may no longer
  be an integral basis.}
modulo $p$, so that this exception does not occur. It suffices to
replace one of the $n$ elements $\xi_i$ by $\xi_i+p$.  We know from
the binomial theorem that
\[(\xi_i+p)^{p^\nu}-(\xi_i+p) \equiv (\xi_i^{p^{\nu}}-\xi_i)-p \equiv
  -p \not \equiv 0 \pmod{P^2},\] because all other terms are divisible
by $p^2$ and so by $P^2$. The easiest way of avoiding the occurrence
of this case, without knowing the prime factors of $p$, is, as is
always possible, to assume that the first element of the fundamental
system a priori equals one, and then instead introduce $\xi_1=1+p$,
which does not change the character of the fundamental system modulo
$p$ at all.\footnote{So Hensel is telling us to always assume
  $\xi_1=1+p$, which guarantees that the $\xi_i$ are no longer an
  integral basis, since $1=\frac{1}{p+1}\xi_1$ is an integer. This is
  not mentioned, however, in the statements of the theorems that
  follow.

  The whole argument is in fact (inadvertently?) local: the $\xi_i$
  are only an integral basis at $p$, and all the divisibility proofs
  below consider only the primes above $p$.} Then this exception can
not occur for any prime factor of $p$ because
\[\xi_1^{p^\nu}-\xi_1=(1+p)^{p^\nu}-(1+p) \equiv -p \pmod{p^2}.\]
In this case $w_{\nu}-w_0$ will not contain any prime divisors
from $p$ more than once. 

The result of these quick observations is summarized in the following
theorem.\footnote{In terms of ideals, the theorem can be stated
  thus. Let $1,\xi_2,\dots,\xi_n$ be an integral basis, and let
  $\xi_1=1+p$. Let $I_\nu$ be the ideal generated by
  $\xi_i^{p^\nu}-\xi_i$, $i=1,2,\dots,n$. Then $I_\nu$ is divisible
  exactly once by each prime ideal of degree dividing $\nu$ that
  appears in the factorization of $p$.

  Note that Hensel does not say (or prove) that there are no divisors
  prime to $p$. Later, he seems to assume that this is the case, but
  it does not hold in general.}

\vspace{.5\baselineskip}

\noindent\hspace{\fill}\parbox{0.8\textwidth}{\hspace*{1em}When $u_1,
  \dots, u_n$ are indeterminates, the linear form
  \[w_{\nu}-w_0= u_1(\xi_1^{p^{\nu}}-\xi_1)+\dots
    +u_n(\xi_n^{p^{\nu}}-\xi_n),\] which is of degree $p^\nu$ with
  respect to the elements of the fundamental system
  $\xi_1, \dots, \xi_n$, is divisible by the product of all distinct
  prime divisors of $p$ whose degree $\kappa$ is an exact divisor of
  $\nu$, and contains each of these exactly once.}

\vspace{.5\baselineskip}

From this theorem \marginpar{[141]} we can draw an interesting
conclusion,\footnote{Hensel will use the theorem to derive a criterion
  for $p$ to be ramified in $(\mathfrak G)$. The point is that $p$ is
  unramified if and only if it divides one of the ideals $I_\nu$. Note
  that he knows that this is equivalent to $p$ not dividing the field
  discriminant.} which is of importance for a subsequent\footnote{As
  Petri notes in \cite[2.4]{Petri}, it is unclear which ``subsequent
  investigation'' Hensel has in mind.} investigation:

If the prime number $p$  contains prime factors which are pairwise
distinct, so that 
\[p \sim P_1P_2 \dots P_h,\]
then we can always find a linear form $w_{\nu}-w_0$ divisible by the
prime $p$, or, equivalently, such that the $n$ congruences: 
\[\xi_i^{p^{\nu}} \equiv \xi_i \pmod{p} \qquad\qquad (i=1,2, \dots, n)\]
are all satisfied. Given the theorem, we only need to choose $\nu$ to
be the least common multiple of the $h$ degrees $\kappa_1, \dots ,
\kappa_h$ of the prime divisors $P_1, \dots , P_h$. 

On the other hand, if $p$ contains a prime divisor more than once,
then none of the linear factors\footnote{Sic, but he wants to say
  ``linear forms.''} $w_{\nu}-w_0$ is divisible by $p$,
because none of those differences can contain a multiple factor of $p$
more than once.  Thus we have the
following theorem: 

\vspace{.5\baselineskip}

\noindent\hspace{\fill}\parbox{0.8\textwidth}{\hspace*{1em}The prime
  $p$ decomposes in the domain $(\mathfrak{G})$ into a product of
  distinct prime divisors if and only if at least one of the linear
  forms $w_{\nu}-w_0$ is divisible by $p$.}

\vspace{.5\baselineskip}
  
Since the number $p$ is a divisor of the domain 
discriminant\footnote{``Gattungsdiscriminante.''} when and only when
it contains at least one multiple prime factor, we can state as a
corollary of the previous result the following theorem:\footnote{We
  have translated the theorem as Hensel states it, but the statement
  seems incorrect. What he had proved is that $p$ is unramified if and
  only if it divides one of the forms $w_\nu-w$ (equivalently, one of
  the ideals $L_\nu$). The negation would then say that $p$ is
  ramified if and only if it divides none of them. It is also unclear
  why he brings in ``a fractional power of $p$.''}

\vspace{.5\baselineskip}

\noindent\hspace{\fill}\parbox{0.8\textwidth}{\hspace*{1em}The prime
  $p$ is contained in the discriminant of the domain $(\mathfrak{G})$
  if and only if one of the linear forms $(w_{\nu}-w_0)$ is divisible
  by a fractional power of $p$, but not divisible by $p$
  itself.}

\vspace{.5\baselineskip}

Now let $\kappa$ be an arbitrary whole number. We will form the product:
\[F_\kappa(w_0)=
  \frac{(w_\kappa-w_0) \prod(w_{\frac{\kappa}{qq'}}-w_0)
    \prod(w_{\frac{\kappa}{qq'q''q'''}}-w_0) \dots}
  {\prod(w_{\frac{\kappa}{q}}-w_0)\prod(w_{\frac{\kappa}{qq'q''}}-w_0)
    \dots},
\]
where $q$, $q'$, $q''$, \dots are the distinct prime factors of
$\kappa$. More simply,
\[F_\kappa(w_0)= \prod_{d|\kappa}(w_d-w_0)^{\varepsilon_d},\] where
$\varepsilon_d= \pm 1$, according to whether the divisor of $\kappa$
complementary to $d$ is a product of an \marginpar{[142]} even or odd
number of \emph{distinct} prime factors $q$, $q'$, $q''$ \dots{} of
$\kappa$, and where $\epsilon_{\delta}=0$ when the ratio
$\frac{\kappa}{d}$ contains repeated prime factors.\footnote{In modern
  terms $\varepsilon_d=\mu(\kappa/d)$, where $\mu$ is the M\"obius
  function. It apparently was introduced by M\"obius in 1832, but was
  clearly not part of the standard toolkit.}  This quotient is a
rational function\footnote{In terms of ideals, the $F_\kappa(w_0)$
  correspond to fractional ideals
  $F_\kappa = \prod_{d|\kappa} I_d^{\epsilon_d}$, where the $I_d$ are
  as above.} of the elements $(\xi_1, \dots, \xi_n )$ of the
fundamental system of $(\mathfrak{G})$. Its dimension\footnote{Hensel
  means the degree of the rational function in the symbols $\xi_i$.}
with respect to these elements is
\[g(\kappa)=p^\kappa- \sum_q p^{\frac{\kappa}{q}}+\sum_q
  p^{\frac{\kappa}{qq'}}- \dots = \sum_{d|\kappa}
  \varepsilon_{d}p^d. \] We immediately recognize that it [namely,
$F_\kappa(w_0)$] is equivalent to\footnote{The claim is that this
  rational function is equivalent, in the sense of Kronecker, to a
  product of prime divisors. In Dedekind's terms, Hensel is saying
  that the ideal $F_\kappa$ is the product of these primes. This is
  incorrect, since Hensel ignores completely the primes that do not
  divide $p$. In other words, the argument continues to be local.}
the product of all distinct prime divisors of $p$ whose degree is
exactly equal to $\kappa$.  If\footnote{First, divisors of degree not
  dividing $\kappa$ do not divide any of the forms that make up
  $F_\kappa(w_0)$.}  $\bar{P}$ is a prime divisor of $p$ whose degree
$\bar{\kappa}$ is not a divisor of $\kappa$, then $\bar{P}$ is
contained in neither the numerator or denominator of
$F_\kappa(w_0)$. If\footnote{Next, divisors of degree dividing
  $\kappa$ but unequal to $\kappa$ cancel out. This is the point of
  the complicated quotient.} however, $\bar{\kappa}$ is a divisor of
$\kappa$ then one shows exactly as in the corresponding
question\footnote{The reference is to the formula for the $n$-th
  cyclotomic polynomial $\Phi_n(x)$ in terms of the polynomials
  $x^d-1$ for $d$ dividing $n$. See, for example,
  \cite[p.~285]{KronNT}, where the notation $\varepsilon_d$ is also
  used.}  in the theory of cyclotomic equations, that $\bar{P}$ occurs
in the denominator just as often as it occurs in the numerator of
$F_\kappa(w_0)$. If\footnote{Divisors of degree $\kappa$ occur exactly
  once.}  $\bar{\kappa}=\kappa$ then $\bar{P}$ is contained once and
only once in the numerator of the $F_\kappa(w_0)$, in the linear form
$(w_\kappa-w_0)$, and thus our claim is proved.\footnote{Paragraph
  break inserted here to improve readability. Note that what is
  missing here is any attempt to deal with primes that are not
  divisors of $p$. See the numerical example below.}

If then $P^{(1)}, P^{(2)}, \dots P^{(\lambda_\kappa)}$ are all the
distinct prime factors of $p$ whose degree is equal to $\kappa$, then
$F_\kappa(w_0)$ is equal\footnote{As a divisor in Kronecker's sense.}
to their product. This means, we have an equivalence:
\[F_\kappa(w_0)=\prod_{d|\kappa}(w_d-w_0)^{\epsilon_d} \sim P^{(1)}
  P^{(2)} \dots P^{(\lambda_\kappa)}.\] Taking norms,\footnote{The
  norm of a prime of degree $\kappa$ is of course $p^\kappa$. The
  equation is not actually true, since Hensel is silently ignoring
  primes that are not divisors of $p$.} it
follows\footnote{The equation effectively defines the number
  $L_\kappa$.  Numerical examples show (see below) that the norm need
  not be a power of $p$, so we should take $L_\kappa$ as the $p$-adic
  valuation of the norm instead.} that
\[N(F_\kappa(w))=p^{\kappa\lambda_\kappa}=p^{L_\kappa}\] [So we
get\footnote{I have left the statement of both theorems in Hensel's
  terms, ``dimension'' and ``degree'' unchanged. The first means the
  degree of the rational function, while the second means the residual
  degree of the corresponding divisor. The ``dimension'' is just
  $g(\kappa)$, which we can easily compute in any case.}]

\vspace{.5\baselineskip}

\noindent\hspace{\fill}\parbox{0.8\textwidth}{\hspace*{1em}The form
  $F_\kappa(w)$, which has dimension $g(\kappa)$ with respect to
  $\xi_1, \dots , \xi_n$, has degree $\kappa\lambda_\kappa$, where
  $\lambda_\kappa$ is the number of distinct prime factors of $p$ of
  degree $\kappa$. If no prime factor of degree $\kappa$ exists, then
  $\lambda_\kappa=0$.} 

\vspace{.5\baselineskip}



If we construct the $n$ forms
\[F_\kappa(w_0) \qquad\qquad (\kappa=1,2, \dots, n),\] we know that
their degree\footnote{``Ordnungszahlen.'' That is not quite the right
  word, since only the $p$-part of the norm has been
  computed. $L_\kappa$ is actually the $p$-adic valuation of the
  norm.} $L_\kappa$ equals $\kappa\lambda_\kappa$; from the previously
proven theorem~(D.), the 
prime $p$ is an inessential divisor of all \marginpar{[143]} equation
discriminants $\mathfrak{D}(\xi_0)$ from $(\mathfrak{G})$ when at
least one of the inequalities
\[\lambda_\kappa>\bar{g}(\kappa)=\frac{1}{\kappa}g(\kappa) 
\qquad\text{ or}\qquad \kappa\lambda_\kappa>g(\kappa)\] holds. Now
since $\kappa\lambda_\kappa$ is the degree of the form $F_\kappa(w_0)$
and $g(\kappa)$ is the dimension with respect to $\xi_1, \dots ,
\xi_n$, we can state the previously found result in a more elegant and
simple form:\footnote{A
  numerical example is clarifying. Suppose $K$ is the number field
  obtained by adjoining a root $\alpha$ of the polynomial
  $x^4 +x^3 + 6x^2 + 2x + 12$. (Global field 4.0.13564.1 in
  \cite{LMFDB}.) An integral basis is then
  $(1,\xi_2=\alpha,\xi_3=\frac12(\alpha^2+\alpha^3),\xi_4=\alpha^3)$. Let
  $p=2$ and set $\xi_1=1+p=3$. We want to consider the ideals $I_k$
  generated by $\xi_i^{2^k}-\xi_i$. Then
  $N(F_1)=N(I_1)=24=2^3\cdot 3$, so $L_1=1\lambda_1=3$. We have
  $N(F_2)=N(I_2I_1^{-1})=1$ and $N(F_3)=N(I_3I_1^{-1})=1$, so
  $L_2=L_3=0$. Finally, $N(F_4)=N(I_4I_2^{-1})=11$ is not divisible by
  $2$, so $L_4=0$. This tells us $2$ is divisible by three prime
  ideals of degree $1$ and by none of degrees $2$, $3$, or $4$. Since
  there are only two irreducible polynomials of degree $1$ in
  $\F_2[x]$, it follows that $2$ is a common inessential discriminant
  divisor in this field. (In fact, the factorization is
  $(2)=\mathfrak p_1^2 \mathfrak p_2\mathfrak p_3$ with all prime
  factors of degree $1$.)

  Notice that in this case $2$ is ramified. Hensel claimed above that
  this happens if and only if $2$ does not divide any of the ideals
  $I_k$, which is easy to check is the case. So for this field $2$ it
  is both an essential and an inessential\label{ess_iness} divisor of
  the discriminant.}

\vspace{.5\baselineskip}
 
\noindent\hspace{\fill}\parbox{0.8\textwidth}{\hspace*{1em}The prime
  $p$ is a common inessential divisor of the equation discriminants
  $\mathfrak{D}(\xi_0)$ of the ring of integers $(\mathfrak{G})$ if
  among the forms
 \[F_\kappa(w_0)= \prod_{d|\kappa}(w_d-w_0)^{\epsilon_d}\qquad\qquad
 (\kappa=1, \dots, n)\] at least one exists whose dimension with
 respect to the elements $\xi_1, \dots , \xi_n$ of the fundamental
 system is smaller than the degree, that is than the exponent
 $L_\kappa$ of $p$ in the equation 
 \[N(F_\kappa(w))=p^{L^\kappa}.\]}

\vspace{\baselineskip}

\begin{center}
  \S3
\end{center}

\vspace{\baselineskip}  

The question of common inessential discriminant divisors can now be
handled in an entirely different fashion, leading to an entirely
different criterion for them to occur.

Let
\[\xi_1^{(0)}, \dots, \xi_n^{(0)}\]
be a fundamental system\footnote{I.e., an integral basis. The
  subscripts $(0)$ have been added because Hensel is about to consider
  conjugates.} for the field $(\mathfrak{G})$ and let
\begin{align*}
w^{(0)} &=u_1\xi_1^{(0)}+ \dots+u_n\xi_n^{(0)},\\
w^{(1)} &=u_1\xi_1^{(1)}+ \dots+u_n\xi_n^{(1)},\\
\vdots & \\
w^{(n-1)}& =u_1\xi_1^{(n-1)}+ \dots+u_n\xi_n^{(n-1)}
\end{align*}
be the fundamental forms for the field $(\mathfrak{G})$ and its
conjugates. Then the discriminant of the fundamental equation is
\[D=\prod_{a \neq \beta}(w^{(a)}-w^{(\beta)})\qquad\qquad (a,\beta =
  0,1, \dots , n-1).\] This is a homogeneous function of
$u_1, \dots, u_n$ with integer coefficients. The greatest common
divisor of all of these coefficients is a whole number,
\marginpar{[144]} which I have shown in the previous work\footnote{As
  usual, this is \cite{Hensel1894a}.} (page 78) agrees with the field
discriminant, which is to say, with the square of the determinant
\[|\xi_i^{(\kappa)}|^2\qquad\qquad
\genfrac{}{}{0pt}{}{(i=1,2, \dots, n)}{\qquad(\kappa=0,1, \dots,
  n-1)}.\] 
So we have\footnote{There are no equations in this section marked (1.)
  or (2.).}
\[\tag{3.} D(u_1, \dots, u_n)=\Delta^2(u_1, \dots, u_n)
  |\xi_i^{(\kappa)}|^2,\] where $\Delta(u_1, \dots, u_n)$ is a
homogeneous function\footnote{In fact $\Delta$ is the ``index form,''
  i.e., in computes the index of $\Z[\xi^{(0)}]$ in the ring of
  integers. For his example of a cubic field in which $2$ is a common
  inessential discriminant divisor in \cite{Ded1878}, Dedekind
  computed it explicitly.} of $u_1, \dots u_n$, whose dimension is
clearly equal to $\frac{n(n-1)}{2}$, and whose coefficients no longer
have any common divisors, which is to say $\Delta(u_1,\dots, u_n)$ is
a primitive polynomial function in $u_1, \dots, u_n$.

Now suppose that instead of $w^{(0)}$ we choose a number from the domain
\[\xi^{(0)}=a_1\xi_1^{(0)}+\dots +a_n\xi_n^{(0)}\] together with its
conjugates. Then we obtain its equation discriminant if in (3.) we
substitute the indeterminates $u_1, \dots, u_n$ by the whole numbers
$a_1, \dots, a_n$. The prime $p$ is an inessential divisor of the
discriminant if and only if it is contained in $\Delta(a_1,\dots,
a_n)$. So we have the following theorem:\footnote{This is known as
  Hensel's criterion for common inessential discriminant divisors. The
  idea is to compute the index form and then check that its values are
  always divisible by $p$. Hensel will state it first, then give an
  explicit way to test a form to see if all its values are indeed
  divisible by $p$, then summarize the whole thing into a theorem.}

\vspace{.5\baselineskip}
 
\noindent\hspace{\fill}\parbox{0.8\textwidth}{\hspace*{1em}The prime
  $p$ is a common inessential equation discriminant divisor for the
  field $(\mathfrak{G})$ if and only if the primitive form
  \[ \Delta(u_1,\dots, u_n)\]
  is divisible by $p$ for all integer values of the indeterminates
  $u_1, \dots, u_n$.}

\vspace{.5\baselineskip}
 
The question of when a polynomial form has values divisible by $p$ for
all integer values of the indeterminates is fully answered by the
following theorem: 

\vspace{.5\baselineskip}
 
\noindent\hspace{\fill}\parbox{0.8\textwidth}{\hspace*{1em}A form
  $U(u_1, \dots, u_n)$ has value divisible by a prime $p$ for all
  integer values of the indeterminates if it contains the
  module system:
  \[\tag{4.} (p; u_1^p-u_1,\dots, u_n^p-u_n);\]
  that is, when $U$ can be written in the form 
  \[\tag{4$^{\mathrm a}$.} U(u_1, \dots,
    u_n)=pU_0+(u_1^p-u_1)U_1+\dots +(u_n^p-u_n)U_n,\] 
  where $U_0, U_1, \dots, U_n$ are integral polynomials in
  $u_1, \dots, u_n$.}

\vspace{.5\baselineskip}
 
\marginpar{[145]} This theorem can most easily be proved\footnote{The
  proof is in this and the next paragraph.} through
induction.\footnote{In fact the proof is constructive: take the form,
  divide it by $u_1^p-u_1$, look at the coefficients of the resulting
  polynomial in $u_1$, rinse, repeat.} It is obviously true when no
variable is present; we now assume that it is proved for the case of
$n-1$ variables $(u_2,\dots, u_n)$ and prove it for $n$ variables. If
the form $U(u_1,\dots, u_n)$ has degree higher that $p-1$ in $u_1$,
then it can be reduced modulo $u_1^p-u_1$ to another form
$\bar{U}(u_1, \dots, u_n)$, whose degree in $u_1$ is at most equal to
$p-1$, since clearly
\[u_1^p \equiv u_1, u_1^{p+1} \equiv u_1^2, \dots, u_1^{p+i} \equiv
u_1^{i+1} \pmod{(u_1^p-u_1)}.\] 
If we write the functions according to the powers of $u_1$, we see get
congruence: 
\[\tag{5.} U(u_1, \dots, u_n) \equiv \bar{U} =
\bar{U}_0u_1^{p-1}+\bar{U}_1u_1^{p-2}+\dots +\bar{U}_{p-1}
\pmod{(u_1^p-u_1)},\] and the function $\bar{U}(u_1, \dots, u_n)$ will
be divisible by the prime $p$ for every whole integer value of $(u_1,
\dots, u_n)$ when the same is also the case for $U(u_1, \dots, u_n)$
and conversely, since they differ by a multiple of $u_1^p-u_1$, which
for every whole value of $u_1$ is itself a multiple of $p$ according
to \emph{Fermat}'s theorem.

Now if we give $u_2, \dots, u_n$ integer values $a_2, \dots, a_n$,
$\bar{U}_0, \dots \bar{U}_{p-1}$ become equal to integers
$A_0, \dots , A_{p-1}$. The resulting expression for $\bar{U}$ 
\[ A_ou_1^{p-1} + A_1u_1^{p-2} + \dots + A_{p-1},\]
must be divisible by $p$ when we let $u_1$ be equal to each of the $p$
incongruent numbers $0,1, \dots , p-1$. But an expression of degree
$p-1$ can only vanish modulo $p$ for $p$ incongruent values of $u_1$
if \emph{all} its coefficients\footnote{I.e., a polynomial of degree
  $p-1$ cannot have $p$ roots.} are divisible by $p$. So it follows
that $\bar{U}(u_1, \dots, u_n)$ is only divisible by $p$ for all
integer value systems if the same is true for the $p$ coefficients
$\bar{U}_0(u_2, \dots, u_n), \dots, \bar{U}_{p-1}(u_2, \dots, u_n)$,
which are functions only of $u_2, \dots, u_n$. If this is the case,
then according to inductive assumption all of these coefficients
contain  the divisor system $p; u_2^p-u_2, \dots , u_n^p-u_n$. The
same now follows for the whole form $\bar{U}(u_1, u_2, \dots u_n)$,
and from equation (5) it follows that the function being investigated
contains the divisor system
\[(p; u_1^p-u_1, \dots, u_n^p-u_n),\] 
\marginpar{[146]} since it differs from the previous by only a
multiple of $(u_1^p-u_1)$; and thus the theorem is proved.

With the help of this theorem we now have the following criterion for
the occurrence of a common inessential discriminant divisor:

\vspace{.5\baselineskip}
 
\noindent\hspace{\fill}\parbox{0.8\textwidth}{\hspace*{1em}If
\[\Delta^2(u_1, \dots, u_n )\]
is the discriminant of the fundamental equation of the
domain $(\mathfrak{G})$ with its numerical factor removed, then the
prime $p$ is a common inessential 
divisor of the equation discriminants for $(\mathfrak G)$ if and only
if the primitive polynomial form $\Delta(u_1, \dots, u_n)$
of degree $\frac{1}{2}n(n-1)$ contains the divisor system 
\[(p; u_1^p-u_1, \dots, u_n^p-u_n),\] (in the sense of
\emph{Kronecker}'s Festschrift), that is, if there is an equation
\[ \Delta(u_1, \dots, u_n)= U_0p+U_1(u_1^p-u_1)+ \dots +U_n(u_n^p-u_n),\]
where $U_0,U_1, \dots, U_n$ are polynomials in $u_1, \dots ,
u_n$.}

\vspace{.5\baselineskip}

In his Festschrift for \emph{E. E. Kummer}'s doctorate anniversary,
\emph{L. Kronecker} mentioned\footnote{In
  \cite[\S25]{Grundzuge}. Kronecker's \emph{Gundz\"uge} was first
  published as to celebrate the fiftieth anniversary of Kummer's
  doctorate.} the possibility of such common inessential discriminant
divisors, adding that this occurs, for example, when the primitive
form $\Delta(u_1, \dots, u_n)$ can be expressed as a homogeneous
polynomial function of $(u_i^p-u_i)$.  The preceding simple
observations show us\footnote{This paragraph is difficult to
  understand, but we think we caught the basic meaning. Hensel is saying
  that Kronecker observed that this condition was sufficient but that
  what he has added is that it is also necessary, since the content of
  the discriminant form is exactly the field discriminant.} that this
is the essential issue, as long as we add a term of the form
$pU_0(u_1, \dots, u_n)$ to every expression. But to complete the proof
we required the result from the previous paper (Page 78) that the form
$\Delta(u_1, \dots, u_n)$ is primitive, which means that the
discriminant $D(u_1 \dots, u_n)$ of the fundamental equation contains
no numerical divisor beyond the field discriminant.

If we know the primitive form $\Delta(u_1, \dots, u_n)$, then it will
be very easy to determine whether or not $p$ is a common inessential
discriminant divisor in $(\mathfrak{G})$. We first reduce the
coefficients of this form to their smallest remainder modulo $p$ and
all exponents of $u_1, \dots , u_n$ larger than $p-1$ to their
smallest remainder modulo $p-1$; \marginpar{[147]} $p$ is a common
inessential divisor of the discriminants of $(\mathfrak{G})$ if and
only if the resulting form is identical to $0$.

To illustrate this way of handling the problem, I will now give the
following simple example,\footnote{In modern language, Hensel looks at
  the cyclotomic field corresponding to a prime number
  $\nu=3\mu+1$. This has a unique cyclic cubic subfield
  $\mathfrak G_3(\epsilon)$ generated by
  Gaussian periods $\epsilon_1,\epsilon_2,\epsilon_3$. In his thesis,
  Hensel found a sufficient condition for such fields (and
  more general versions of them) to have common inessential
  discriminant divisors.} which was considered from a different point
of view in my doctoral dissertation.

\vspace{.5\baselineskip}
 
\noindent\hspace{\fill}\parbox{0.8\textwidth}{\hspace*{1em}Let $\nu$
  be an arbitrary real prime of the form $3\mu+1$ and
  $\mathfrak{G}_3(\epsilon)$ be the field generated by the three
  $\mu$-fold periods 
  $\epsilon_1, \epsilon_2, \epsilon_3$ of the $\nu$th roots of
  unity. We want to find the common inessential
  discriminant divisors of this field.} 

\vspace{.5\baselineskip}
 
The three periods $\epsilon_1, \epsilon_2, \epsilon_3$ form a
fundamental system for the field\footnote{Hensel doesn't say, but he
  knows and uses, that the three periods are cyclically permuted by
  the Galois action. The three forms $w_1$, $w_2$, $w_3$ below are,
  of course, Galois conjugates.} $\mathfrak{G}_3(\epsilon)$. Let 
\begin{align*}
  w_1 &=u_1\epsilon_1+u_2\epsilon_2+u_3\epsilon_3,\\
  w_2 &=u_1\epsilon_2+u_2\epsilon_3+u_3\epsilon_1,\\
  w_3 &=u_1\epsilon_3+u_2\epsilon_1+u_3\epsilon_2,
\end{align*}
so that the product $(w_1-w_2)(w_2-w_3)(w_3-w_1)$ will be the square
root of the discriminant of the fundamental equation. By using the
known expressions for the resolvent of the cubic period equation, we
get without difficulty the expression
\[\prod (w_i-w_{i+1})=-\nu(\alpha\Delta_1+\beta\Delta_2).\]
Here $\alpha$ and $\beta$ are the integers which occur in the
decomposition of $\nu$ into its prime factors in the field of the
third roots of unity, so that 
\[\tag{6.} \nu = (\alpha+3\beta\rho)(\alpha+3\beta\rho^2) =
\alpha^2-3\alpha\beta+9\beta^2\qquad\qquad (\rho^2+\rho+1=0),\]  
and $\Delta_1$ and $\Delta_2$ are the primitive forms:
\begin{align*}
  \Delta_1 &=(u_1-u_2)(u_2-u_3)(u_3-u_1)= \sum_{i=1}^3u_i u_{i+1}^2 -
  \sum_{i=1}^3u_i^2 u_{i+1}, \\
  \Delta_2 &= \sum(u_i^3-3u_i u_{i+1}^2)+6u_1u_2u_3.
\end{align*}
Since the form $\Delta_1$ is only of degree $2$ with respect to the
quantities $u_1, u_2, u_3$, only the number $2$ can occur as a common
discriminant divisor. If we reduce the forms $\Delta_1$ and
$\Delta_2$ modulo the divisor system
\[M=(2, u_1^2-u_1, u_2^2-u_2, u_3^2-u_3),\]
we quickly see that $\Delta_1$ contains the same and we get the
following congruence \marginpar{[148]} for the primitive form
$(\alpha\Delta_1+\beta\Delta_2)$:
\[\alpha\Delta_1+\beta\Delta_2 \equiv
\beta[(u_1+u_2+u_3)+(u_1u_2+u_2u_3+u_3u_1)] \pmod{M},\] 
which means the number $2$ is only a common discriminant divisor when 
\[\beta \equiv 0 \pmod{2}.\]

We can now give this result a more elegant form. The decomposition of
$\nu$ in (6.) allows us to represent $4\nu$ as:
\[4\nu=(2\alpha-3\beta)^2+27\beta^2=A^2+27B^2,\]
so that the number $4\nu$ can always be expressed as $A^2+27B^2$, and
the uniqueness of the decomposition in (6.) shows 
us that our depiction is also unique. From the two equations  
\[A=2\alpha-3\beta, B=\beta\]
it follows that $A$ and $B$ are only divisible by $2$ if the same is
true for $\beta$, which means that in this case not only $4\nu$ but
also $\nu$ itself can be expressed in the form $A^2+27B^2$. So we have
the following theorem:\footnote{Note that this is stated incorrectly
  in \cite[2.2.1, item 3]{Nark}.} 

\vspace{.5\baselineskip}
 
\noindent\hspace{\fill}\parbox{0.8\textwidth}{\hspace*{1em}If
  $\nu=3\mu+1$ is an arbitrary real prime, then the number $2$ is 
  a common inessential discriminant divisor in the field
  $\mathfrak{G}_3(\epsilon)$, if and only if the number $\nu$ can
  be written in the form  
  \[\nu=A^2+27B^2.\]}

\vspace{.5\baselineskip}

For primes less than two hundred whose fields of roots of unity
contain cubic subfields,\footnote{Hensel's actual sentence is
  something like: ``For primes $\nu$ in the first and second hundred
  this case happens for cubic period equations which are formed by
  roots of unity of order\dots''}
  this case occurs for
\[31, 43, 109, 127, 157, 189.\]

\vspace{2\baselineskip}

\begin{center}
  \S4
\end{center}

\vspace{\baselineskip}  

In his discussion mentioned above of the common inessential
discriminant divisors of a domain, \emph{Kronecker} calls attention to
the remarkable circumstance that this can be
eliminated,\footnote{\label{KL}Hensel is unclear on what exactly can
  be ``eliminated'' here. Suppose $p$ is a common inessential
  discriminant divisor for a number field. By
  the previous theorem, the index form is a primitive form in the
  variables $u_i$ which becomes divisible by $p$ whenever we replace
  each $u_i$ by an integer $a_i\in\Z$. Kronecker's observation is that
  we can obtain a value that is not divisible by $p$ if we allow the
  values $a_i$ to be algebraic integers in a larger field (which
  Hensel calls $\Gamma(\zeta)$). (In terms of theorem~A, we are
  replacing $\F_p$ by a finite field extension to get more irreducible
  polynomials.) 

  Hensel proposes to give a proof of Kronecker's remark and to explain
  how to find a field $\Gamma(\zeta)$. In fact, he will claim that he
  can take $\Gamma(\zeta)$ to be a subfield of a cyclotomic field
  $\Q(\mu_\ell)$ with $\ell$ a prime.

  Denote the original number field by $K$. In modern terms, Hensel
  wants to construct an auxiliary number field $L$, contained in a
  prime-order cyclotomic field. If we assume that $K$ and $L$ are
  linearly disjoint, then the integral basis
  $\{\xi_1,\xi_2,\dots,\xi_n\}$ of $K$ over $\Q$ will also be an
  integral basis for $KL$ over $L$ and the relative discriminant
  $d(KL/L)$ is the ideal in $\mathcal{O}_L$ generated by
  $d_K$. Hensel's result then means that $p$ is no longer a common
  index divisor for $KL/L$. That is, there exists an element
  $\theta\in KL$ such that the index of $\mathcal{O}_L[\theta]$ in
  $\mathcal{O}_{KL}$ is not divisible by $p$.} if the coefficients
$u_1, u_2, \dots u_n$ of the linear form
\[ w_0 = u_1\xi_1 + \dots + u_n\xi_n \] are no longer in the domain of
the real integers, but rather in the bigger \marginpar{[149]} realm of
the algebraic numbers from another domain
$\Gamma (\zeta)$.\footnote{This observation is in
  \cite[\S25]{Grundzuge} (p.~384 in volume II of \cite{KronWerke}).}
\emph{Kronecker} does not prove this interesting theorem, however, nor
does he specify how the ring of integers should be chosen such that we
can avoid the occurrence of a common discriminant
divisor. \emph{Kronecker} did not return to this subject later, and so
far I cannot find\label{cannotfind} any hint of a proof in his
papers.\footnote{Did Hensel wait until after Kronecker's death to
  publish these results because he expected to find such a proof?}

I would like to briefly touch on this point to show how to choose the
field $\Gamma(\zeta)$ for the coefficients $u_1, \dots , u_n$ so that
the prime $p$ is not a common discriminant divisor, and to determine
the field of smallest degree for the adjunction. Finally, we wish to
show that for this goal we only need the simplest algebraic numbers,
which stem from the roots of unity of prime degree.

The following theorem, which is a simple extension of the one
established in the previous sections, leads to these results:

\vspace{.5\baselineskip}

\noindent\hspace{\fill}\parbox{0.8\textwidth}{\hspace*{1em}Let 
\[\tag{1.} F_1(u_1), F_2(u_2), \dots F_n(u_n)\]
be $n$ integral polynomials, each in one variable $u_1, \dots, u_n$,
such that each $F_i(u_i)$ modulo $p$ has as many incongruent integer
roots $z_i$ as its degree. Further, let $F(u_1, \dots, u_n)$ be an
integral polynomial in all $n$ indeterminates $u_1, \dots, u_n$. Then
the congruence
 \[\tag{2.} F(z_1, z_2, \dots, z_n)
 \equiv 0 \pmod{p}\]  
 holds for all congruence roots $z_1, \dots, z_n$
 of the $n$ functions (1.) if and only if $F(u_1, \dots, u_n)$
 contains the divisor system
 \[(p; F_1(u_1), \dots, F_n(u_n))\]
 in the sense of \emph{Kronecker}'s theory.
}

\vspace{.5\baselineskip}

The proof of this theorem can be carried out in the same way as
before, since the argument was based solely\footnote{Working over
  $\F_p$, the argument boils down to the observation that modulo
  $F_1(u_1)$ the polynomial $F$ becomes a polynomial of degree
  $\deg(F_1)-1$ in $u_1$ whose coefficients are polynomials in the
  other variables. Replace $u_2,\dots, u_n$ with arbitrarily chosen
  roots $z_2,\dots,z_n$. We get a polynomial in $u_1$ which is zero
  for all possible choices of $z_1$. Since the number of choices is
  $\deg(F_1)$, so higher than the degree of the specialized
  polynomial, this polynomial must be identically $0$. Hence each of
  the coefficient polynomials has the property that it is zero for all
  choices of $z_2,\dots,z_n$. Now use induction.} on the fact that the
congruences of degree $p$
\[F_i(u_i)=u_i^p-u_i \equiv 0 \pmod{p}\]
have exactly $p$ incongruent roots modulo $p$, so as many as their degree,
together with the fact that $p$ is a prime number.

\marginpar{[150]} This theorem can finally also be
extended\footnote{Now we allow the polynomials to have coefficients in
  some ring of integers and replace $p$ by one of its prime divisors;
  this amounts to working over a finite extension of $\F_p$, and the
  argument goes through as before.} in the
following way: we can assume that the coefficients of the $n+1$
functions $F_1(u_1)$, \dots , $F_n(u_n)$, $F(u_1, \dots, u_n)$ are no
longer real integers, but rather elements of the ring of integers of a
field $\Gamma(\zeta)$ determined by an arbitrary algebraic integer
$\zeta$; only now we must work modulo a prime divisor $p(\zeta)$ of
$p$ instead of the element $p$, which inside $\Gamma(\zeta)$ 
loses the property of being indecomposable.  Now if the
functions\footnote{Hensel does not just assume the coefficients are
  now in $\Z[\zeta]$, but rather indicates this explicitly in his
  notation.} 
\[F_i(u_i, \zeta)\] 
modulo $p(\zeta)$ contain the same number of incongruent roots
$\zeta_i$ inside the domain $\Gamma(\zeta)$ as their degree indicates,
then we can show as before that for a polynomial
$F(u_1, \dots, u_n,\zeta)$, the congruences
\[F(\zeta_1, \zeta_2, \dots, \zeta_n; \zeta) \equiv 0
\pmod{p(\zeta)}\]
will hold for all value systems $(\zeta_1,\dots, \zeta_n)$ if and only
if the polynomial $F(u_1, \dots, u_n, \zeta)$ can be represented by
the elements of the divisor system 
\[(p(\zeta); F_1(u_1, \zeta), \dots, F_n(u_n, \zeta))\] in a
homogeneous and linear way with integer\footnote{Hensel probably means
  that the coefficients are polynomials in the $u_i$ with integer
  coefficients.} coefficients.

We will now make a special assumption, that the \emph{coefficients} of
each of the $n+1$ function $F_i(u_i,\zeta)$ and
$F(u_1, \dots, u_n, \zeta)$ are, as before, real
integers.\footnote{Hensel will assume the polynomials in question have
  rational integer coefficients, but still wants to allow the
  variables to take values in $\Gamma(\zeta)$. He claims that in this
  case the divisor $p(\zeta)$ above can in fact be replaced by $p$.}
To call attention to this assumption, we will denote them as before by
$F_i( u_i)$ and $F(u_1, \dots, u_n)$. The congruence roots $\zeta_i$
of the function $F_i(u_i)$ modulo $p(\zeta)$, however, are still
assumed to belong to the ring of integers $\Gamma(\zeta)$. Now if we
reduce the integral function $F(u_1, \dots, u_n)$ to the smallest
remainder modulo the integral module system
\[(F_1(u_1), \dots,  F_n(u_n) ),\]
we get an integral function $\bar{F}(u_1, \dots, u_n)$ of
$u_1, \dots , u_n$ with \emph{real} integer coefficients, whose degree
in $u_i$ will always be smaller than the degree of 
the function $F_i( u_i)$. This reduced function can only be
represented in a homogeneous way by the elements of the system 
\[(p(\zeta), F_1(u_1), \dots, F_n (u_n)),\]
if all of their coefficients are divisible by $p(\zeta)$, that is to
say by $p$ itself. Thus $F(u_1, \dots u_n)$ contains the divisor
system 
\[(p(\zeta), F_1(u_1), \dots, F_n(u_n))\] if and only
\marginpar{[151]} if it can be represented in a homogeneous and linear
way by the real\footnote{Recall that Hensel uses ``real'' to mean
  ``rational.''} system
\[(p, F_1(u_1), \dots, F_n (u_n))\]

\vspace{.5\baselineskip}

\noindent\hspace{\fill}\parbox{0.8\textwidth}{\hspace*{1em}The
  function $F(u_1, \dots, u_n )$ vanishes modulo $p(\zeta)$ for all
  congruence roots of the $n$ functions $F_i(u_i)$ if and only if the
  congruence  
  \[\tag{3.} F(u_1, \dots, u_n) \equiv 0~~~~[\mathrm{modd}{(p, F_1(u_1), \dots,
      F_n(u_n))}]\] is satisfied.}

\vspace{.5\baselineskip}

We can use this theorem to solve easily the question posed at the
beginning of this paper. Take, as in the first paragraph of this
paper,  
\[w_0=u_1\xi_1+u_2\xi_2+ \dots+u_n\xi_n\]
the fundamental form of the ring of integers $(\mathfrak{G})$, and let 
$w_1, w_2, \dots , w_{n-1}$ denote the $n-1$ fundamental forms conjugate
to $w_0$. Finally, let
\[\Delta^2(u_1, \dots, u_n)\] be the discriminant of the
fundamental equation freed from its numerical divisors (the field
discriminant\footnote{As before $\Delta$ is the index form.}). Now if
$\Gamma(\zeta)$ is another arbitrary field domain and $p(\zeta)$ is a
prime divisor of the real prime $p$ in $\Gamma(\zeta)$, we can
investigate under which conditions $p(\zeta)$ is a common inessential
divisor of all discriminants $\prod(w_i-w_n)$, by now letting
$u_1, \dots , u_n$ be arbitrary algebraic integers of the domain
$\Gamma(\zeta)$ instead of arbitrary real integers. Equivalently, we
can investigate the conditions under which the primitive form
$\Delta(u_1, \dots, u_n)$ is always divisible by $p(\zeta)$, if we
replace $u_1, \dots , u_n$ by arbitrary integers belonging to the
domain $\Gamma(\zeta)$.

Now if $k$ is the degree of the prime divisor $p(\zeta)$ for the
domain $\Gamma$, then the number of integers\footnote{Hensel always
  thinks in terms of representatives rather than congruence classes;
  so it's ``the number of incongruent integers'' rather than ``the
  number of congruence classes.''} in $\Gamma$ that are incongruent
modulo $p(\zeta)$ is equal to $p^k$. So every number $\zeta$ of
this domain satisfies the congruence:
\[w^{p^k}-u \equiv 0 \pmod{p(\zeta)},\]
and so this congruence contains the same number of
incongruent roots inside of $\Gamma$ as its degree displays. The above
question can now be stated as follows: Under what conditions is the
primitive form $\Delta(u_1, \dots, u_n)$ divisible by $p(\zeta)$ for
the congruence roots modulo $p(\zeta)$ of the $n$ functions  
\[u_1^{p^k}-u_1, u_2^{p^k}-u_2, \dots, u_n^{p^k}-u_n.\] \marginpar{[152]}
This question is directly answered by the last theorem, if we
replace the $n$ functions $F_i(u_i)$ by $u_i^{p^k}-u_i$. So we get the
theorem:\footnote{The statement is confusing because it refers to the
  ``equation discriminants of $(\mathfrak{G})$'' being divisible by a
  divisor in $\Gamma(\zeta)$. See footnote~\ref{KL} for our
  interpretation.} 

\vspace{.5\baselineskip}

\noindent\hspace{\fill}\parbox{0.8\textwidth}{\hspace*{1em}The prime
  divisor $p(\zeta)$ in the field of rationality $\Gamma(\zeta)$ is a
  common inessential divisor of all equation discriminants of
  $(\mathfrak{G})$ if and only if $\Delta(u_1, \dots, u_n)$, the 
  discriminant of $(\mathfrak{G})$ freed of its numerical factor,
  contains the   divisor system 
  \[\tag{4.} P_k=(p; u_1^{p^k}-u_1, \dots, u_n^{p^k}-u_n)\]
  in \emph{Kronecker}'s sense, where $k$ is the degree of $p(\zeta)$ for
  the domain $\Gamma$.}

\vspace{.5\baselineskip}

Since the above criterion is solely dependent on the degree $k$ of
$p(\zeta)$, it applies equally to all divisors of $p$ in $\Gamma$ with
the given degree.

It follows from this that in the rationality domain $\Gamma(\zeta)$ the
prime divisor $p(\zeta)$ is not a common inessential discriminant
divisor of $(\mathfrak{G})$ if and only if the following condition is
satisfied: 
\[\Delta (u_1, \dots, u_n) \not \equiv 0 \pmod{(p; u_1^{p^k}-u_1,
  \dots, u_n^{p^k})}.\] 

This result can now be used to decide what assumptions need to be made
on the values taken by coefficients $u_1, \dots , u_n$ of the
fundamental form of $(\mathfrak{G})$, $w_0=u_1\xi_1+ \dots +u_n\xi_n$,
inside a domain 
$\Gamma(\zeta)$, so that the discriminant $D(w_0)$ of the $n$
conjugate values $w_0, w_1, \dots , w_{n-1}$ does not contain \emph{any}
prime factor of $p$ other than those in the discriminant, or equivalently,
whether it is possible to choose values for the unknowns
$u_1, \dots , u_n$ in $\Gamma(\zeta)$ so
that the form $\Delta(u_1, \dots, u_n)$ is relatively prime to $p$.

Clearly we will first need to choose $u_1, \dots , u_n$, so that the
primitive form $\Delta(u_1, \dots u_n)$ is coprime to every prime
divisor of $p$ in $\Gamma(\zeta)$, i.e., so that none of the prime
divisors [of $p$ in $\Gamma(\zeta)$] are common inessential divisors
of the discriminants $D(w_0)$ from $(\mathfrak{G})$. If
\[p=p_1^{\epsilon_1}(\zeta) \dots p_l^{\epsilon_l}(\zeta)\]
is the decomposition of $p$ into its prime factors in $\Gamma(\zeta)$,
and if \[k_1, \dots, k_l\]
\marginpar{[153]} are the degrees of the individual distinct prime
divisors, then $p$ can only have the required properties if the
primitive form $\Delta(u_1, \dots, u_n)$ does not contain any
of the $l$ divisor systems 
\[P_{k_i}=(p; u_1^{p^{k_i}}-u_1, \dots, u_n^{p^{k_i}}-u_n)
\qquad\qquad (i=1,2, \dots, l),\]
where of course we only need to investigate those systems for which
the numbers $k_i$ are distinct. If these conditions are satisfied,
then it is easy to see that for the unknowns $u_1, \dots , u_n$, such
integers $\zeta_1, \dots , \zeta_n$ of the domain $\Gamma(\zeta)$
can\footnote{Hensel claims here that if we know we can make each of a
  set of conditions hold separately, we can make them hold
  simultaneously. The ``Chinese remainder theorem'' argument is given
  in the rest of this paragraph.} be chosen so that the number
$\Delta(\zeta_1, \dots, \zeta_n)$ is co-prime to $p$. If each of the
$p_i(\zeta)$ is not a common inessential divisor of the equation
discriminant $D(w_0)$, then for each $i$ we can find $n$ numbers
$\zeta_1^{(i)}, \dots , \zeta_n^{(i)}$ such that  
\[\Delta(\zeta_1^{(i)}, \dots, \zeta_n^{(i)}) \not \equiv 0
\pmod{p_i(\zeta)}.\] 
If we now consider these numbers for each of $l$ prime divisors of
$p$, we can choose other numbers $\zeta_1, \dots,  \zeta_n$ so that
for each $i$ we have 
\[\zeta_1 \equiv \zeta_1^{(i)}, \zeta_2 \equiv \zeta_2^{(i)}, \dots,
\zeta_n \equiv \zeta_n^{(i)} \pmod{p_i(\zeta)}(i=1,2, \dots, l),\] 
and so for each $i$:
\[\Delta(\zeta_1, \dots, \zeta_n) \equiv \Delta(\zeta_1^{(i)}, \dots,
\zeta_n^{(i)}) \not \equiv 0 \pmod{p_i(\zeta)}.\] This means the number
$\Delta(\zeta_1, \dots, \zeta_n)$ is in fact coprime to
$p$.\footnote{Added a paragraph break here.}

We will call the domain $\Gamma(\zeta)$ a \emph{supplementary domain}
for the domain $\mathfrak{G}(\xi)$ with respect to the prime $p$ if we
can choose values in in $\Gamma(\zeta)$ for the unknowns $u_1, \dots ,
u_n$ in the $n$ conjugate fundamental forms $w_0, w_1, \dots ,
w_{n-1}$ so that the discriminant $\prod(w_i-w_k)$ contains the prime
$p$ no more than the domain discriminant of $(\mathfrak{G})$. Then we
can state the necessary and sufficient conditions for $\Gamma(\zeta)$
to be a supplementary domain for $\mathfrak{G}(\xi)$ with respect to
$p$:

\vspace{.5\baselineskip}

\noindent\hspace{\fill}\parbox{0.8\textwidth}{\hspace*{1em}Let
  $p=p_1^{\epsilon_1}p_2^{\epsilon_2}, \dots, p_l^{\epsilon_l}$ be the
  decomposition of the real prime $p$ into prime factors inside the
  domain $\Gamma$, and let $k_1, k_2, \dots , k_\lambda$ be the
  distinct [residual] degrees [of the $p_i$]. Then $\Gamma(\zeta)$ is
  a supplementary domain for $(\mathfrak{G})$ with respect to the prime
  $p$ if and only if the discriminant of the fundamental
  equation of $(\mathfrak{G})$ freed from its numerical factor does not
  contain any of the $\lambda$ divisor systems: 
  \[P_{k_1}, P_{k_2}, \dots , P_{k_\lambda}.\] }

\vspace{\baselineskip}

\marginpar{[154]} If the adjoined domain $\Gamma(\zeta)$ is \emph{Galois},
then the degrees of all prime divisors of $p$ are equal. Setting $k$
to be the common value of all the degrees, we can replace the more
complicated condition above by the simpler condition that the form
$\Delta(u_1, \dots, u_n)$ does not contain a divisor system 
\[P_k=(p; u_1^{p^k}-u_1, \dots, u_n^{p^k}- u_n).\]  

Further, I would also like to remark that in the above theorem all of
the divisor systems $P_k$ can be omitted if the index $k$ is a
multiple of one of the other numbers $k_1, \dots , k_{\lambda}$. If 
\[\Delta(u_1, \dots, u_n) \not \equiv 0 \pmod{P_k},\]
then a fortiori
\[\Delta(u_1, \dots, u_n) \not \equiv 0 \pmod{P_{ak}};\]
because from the congruence
\[u^{p^{ak}}-u \equiv (u^{p^k}-u)^{p^{(a-1)k}} \equiv 0 \pmod{(p; u^{p^k}-u)}\]
it follows that every divisor system $P_{ak}$ is a multiple of $P_k$:
so $\Delta(u_1, \dots, u_n)$ cannot be divisible by $P_{ak}$ if it does
not contain the system $P_k$.  

\vspace{\baselineskip}

\begin{center}
  \S5
\end{center}

\vspace{\baselineskip}  

We should now discuss which is the supplementary domain
$\Gamma(\zeta)$ of lowest degree for a given domain $(\mathfrak{G})$
with relation to an arbitrary prime $p$.

To this end, we investigate discriminant of the fundamental equation
freed from its numerical factor
\[\Delta(u_1, \dots, u_n )\]
as to its divisibility by the divisor systems
\[P_1, P_2, P_3, \dots, \]
where in each case
\[P_k=(p; u_1^{p^k}-u_1, \dots, u_n^{p^k}-u_n).\]
\marginpar{[155]} If the primitive form does not already contain
the first system 
\[P_1=(p; u_1^p-u_1, \dots, u_n^p-u_n),\]
then $p$ is not at all an inessential divisor for $(\mathfrak{G})$,
which means the supplementary domain of lowest degree is that of
natural integers. If however, $\Delta(u_1, \dots, u_n )$ is divisible
by $P_1$, in the sequence $P_1, P_2, \dots$ we must come at last  to
a divisor system 
\[P_k= (p; u_1^{p^k}-u_1, \dots, u_n^{p^k}-u_n),\]
which is no longer contained\footnote{Hensel is stating a lemma: it is not
  possible that for all $k$ the form $\Delta$ is a linear combination
  of the elements in $P_k$. The proof, given in the rest of this
  paragraph, is easy: no form can be a linear combination of
  polynomials of degree bigger than its own degree unless it is zero,
  so if $\Delta\in P_k$ for large enough $k$ it would have to be
  divisible by $p$. But $\Delta$ is primitive. That puts an upper
  bound on $k$.} in $\Delta$. Now choose $\mu$ large enough so that the
power $p^{\mu}$ is bigger than $\frac{n(n-1)}{2}$, which is to say
larger than the dimension\footnote{As before, Hensel seems to use
  ``dimension'' for the degree of a homogeneous form.} of the
primitive form $\Delta(u_1, \dots,u_n)$.  Then this form cannot be
reduced to one of lower degree modulo the system $(u_1^{p^{\mu}}-u_1,
\dots, u_n^{p^{\mu}}-u_n)$, because all exponents of $u_1, \dots, u_n$
will be smaller than $p^{\mu}$. Therefore, the form $\Delta(u_1,
\dots, u_n )$ could only contain the module system
\[(p; u_1^{p^{\mu}}-u_1, \dots, u_n^{p^{\mu}}-u_n),\]
if all its coefficients were divisible by $p$, which contradicts
the assumption that $\Delta(u_1, \dots, u_n )$ is primitive. Thus, the
form can only contain a finite number of divisor systems $P_1, P_2,
\dots$, and the number of systems it contains will be smaller than
$\mu$, where $p^{\mu}$ is the lowest power of $p$, which is bigger than
$\frac{n(n-1)}{2}.$

Now take $P_k$ to be the first module system of the series $P_1, P_2,
\dots, P_{k-1}, P_k$, which is not contained in $\Delta$, so that each of
the previous contain the form. Now if
\[\varphi(\zeta)=\zeta^k + a_1\zeta^{k-1}+ \dots + a_k =0\]
is an equation of degree $k$, whose left side is also irreducible
modulo $p$,\footnoteB{A polynomial of
  degree $k$ that is irreducible modulo $p$ always exists because the
  number \[\bar{g}(k)=\frac{1}{k}\sum_{d|k}\epsilon_{\delta}p^d\] of
  such functions is never equal to zero.} then the domain
$\Gamma(\zeta)$ that it defines will be\footnote{This claim is proved
  in the next two paragraphs. Notice that Hensel doesn't care if the
  supplementary field is linearly disjoint from his original field;
  that only matters if we want to interpret his result in terms of a
  relative extension $LK/L$ as above. So he can use any equation of
  degree $k$ which is irreducible mod $p$.} a supplementary domain for
$\mathfrak{G}(\xi)$ with respect to $p$, and indeed it will be one of
the lowest possible degree.

Indeed, $\Gamma(\zeta)$ is a supplementary domain of
($\mathfrak{G})$. \marginpar{[156]} The function $\varphi(\zeta)$ is
irreducible modulo $p$, so $p$ is itself a prime inside the domain
$\Gamma(\zeta)$ whose [residual] degree is\footnote{Sic, but Hensel
  means $k$. He says ``ihre Ordnung f\"ur denselben ist $p^k$,''
  literally ``the order for the same is $p^k$,'' so perhaps he means
  the number of congruence classes?} $p^k$. The prime $p$ will be a
common inessential divisor in the rationality domain $\Gamma(\zeta)$
of the equation discriminants of $(\mathfrak{G}$ if and only if the
form $\Delta(u_1, \dots u_n)$ contains the divisor system
$P_k=(p; u_i^{p^k}-u_i)$, which contradicts the previous assumption.

Furthermore, if $\Gamma_1(\eta)$ is a different domain whose degree is
smaller than $k$, and $p(\eta)$ is a prime divisor of $p$, then its
degree $k_1$ is at most equal to the degree of $\Gamma_1(\eta)$, and
so is smaller than $k$. Therefore $p(\eta)$ is a common inessential
divisor for $(\mathfrak{G})$ in the domain $\Gamma_1(\eta)$, because
the form $\Delta(u_1, \dots, u_n)$ contains the divisor system
$P_{k_1}= (p; u_i^{p^{k_1}}-u_i)$, whose index is smaller than
$k$.\footnote{Since all he really needs is for the residual degree of
  at least one of the factors of $p$ to be equal to $k$, the smallest
  possible degree for the supplementary field is realized when $p$ is
  inert and $f=k$.} So
we have the following theorem:

\vspace{.5\baselineskip}

\noindent\hspace{\fill}\parbox{0.8\textwidth}{\hspace*{1em}If
  \[P_k=(p; u_1^{p^k}-u_1, \dots, u_n^{p^k}-u_n)\] is the divisor
  system of lowest degree that is \emph{not} contained in the
  primitive form $\Delta(u_1, \dots, u_n)$, then the smallest
  supplementary domain $\Gamma$ of $\mathfrak{G}(\xi)$ for the prime
  $p$ has degree $k$. Such a domain will be defined by every
  polynomial equation of degree $k$, whose left side is irreducible
  modulo $p$.}

\vspace{.5\baselineskip}

From the proof it follows also that we can only obtain a supplementary
domain defined by an equation of degree $k$ if the left side is
irreducible modulo $p$. Further,\footnote{We don't understand this
  paragraph. Hensel seems to be claiming that there is only one field
  $\Gamma(\zeta)$ of this type, but that is not true. The residue
  fields will all be $\F_{p^k}$, of course, and that may be what he is
  referring to. See the next paragraph, where he admits he doesn't
  really mean what he has said.} the domain $\Gamma(\zeta)$ defined by
the above equation, that is, the totality of rational functions of
$\zeta$ contain all functions, which are irreducible equations modulo
$p$ of $k$th degree, and so we see that there is only one domain
$\Gamma(\zeta)$, which is a supplementary domain of lowest degree with
respect to $\mathfrak{G}(\xi)$.

This last statement should be understood as follows. The algebraic
integers in two fields of degree $k$ [defined by equations] which are
irreducible modulo $p$, are pairwise congruent modulo this prime, so
that for the question we are considering one domain can be substituted
for the other. The algebraic character of the supplementary domains
$\Gamma$ can, however, be very different, which raises the question of
which is the algebraically simplest\footnote{Hensel doesn't say what
  he means by ``algebraically simplest,'' of course. It will turn out
  that he can always choose a cyclotomic field.} supplementary domain
for a given domain $(\mathfrak{G})$ with respect to the prime $p$.

\marginpar{[157]} In the previously discussed
place\footnote{Kronecker's \emph{Grunz\"uge} \cite[\S25]{Grundzuge}.}
in his Festschrift, \emph{Kronecker} considers\footnote{Hensel says
  something like ``Kronecker suggests that for this field\dots for
  which $2$\dots (as was first noted\dots), that\dots{} That was too
  much, so we broke it up.} the field defined by the cubic
equation \[\alpha^3-\alpha^2-2\alpha-8=0,\] for which the number $2$
is a common inessential discriminant divisor (as was first noted by
Mr.~\emph{Dedekind} in the cited\footnote{This is \cite{Ded1878}, but
  Dedekind had actually given this example earlier, in
  \cite{anzeige}.} paper).  [Kronecker suggests] that this prime stops
being an inessential divisor if the domain $\Gamma$ of the third roots
of unity is adjoined to the rational domain.\footnote{With the
  integral basis $\{1,\alpha,4/\alpha=\frac12(\alpha^2+\alpha)-1\}$
  given by Dedekind, the index form in this case is
  $\Delta(u_1,u_2,u_3)=2u_2^3 -u_2^2u_3 - u_2u_3^2 - 2u_3^3$, which is
  clearly always even for integer values of the $u_i$. If, however,
  $\zeta$ is a cube root of unity, then $f(0,\zeta,\zeta^2)=1$, as
  Kronecker says.} For this domain $(\mathfrak{G})$, then, the very
simple \emph{cyclotomic} field of the third roots of unity is a
supplementary domain with respect to the inessential divisor $2$.

This suggests an interesting theorem, that for every domain
$(\mathfrak{G})$ and an arbitrary prime $p$ we can find a
supplementary domain of the greatest algebraic simplicity, namely one
constructed from\footnote{It seems likely that when Hensel says
  ``constructed from roots of unity'' he means a subfield of a
  cyclotomic field. But one could also ask for a cyclotomic field
  rather than a subfield. In fact, in what follows Hensel first finds
  the smallest prime-order cyclotomic field that has the desired
  property, then finds the smallest subfield of that field that still
  has the property.} roots of unity of prime degree. In what follows
we prove this theorem and develop a method to find the smallest such
domain.\footnote{To clarify the argument, we exemplify using the cubic
  field above. We use the integral basis $\{1,\alpha,\beta\}$ where
  $\beta=4/\alpha$.}

Let $(\mathfrak{G})$ be a domain of $n$th degree and let
$\Delta(u_1, \dots, u_n)$ be the discriminant of the fundamental
equation of $(\mathfrak{G})$ freed of its numerical
factor.\footnote{The discriminant of the fundamental equation is
  $-503(2u_2^3 -u_2^2u_3 - u_2u_3^2 - 2u_3^3)^2$. While Hensel never
  says it explicitly, to find his $\Delta$ we first divide by the
  field discriminant (here, $-503$) and then take the square root, so
  in this case
  $\Delta(u_1,u_2,u_3)=2u_2^3 -u_2^2u_3 - u_2u_3^2 - 2u_3^3$, as
  mentioned above. That $2$ is a common inessential discriminant
  divisor follows from
  $\Delta(u_1,u_2,u_3)= -u_3(u_2^2-u_2) - u_2(u_3^2-u_3) +
  2(u_2^3+u_3^3-u_2u_3)$.}  Further let
\[P_{k_1}, P_{k_2}, \dots P_{k_l}\]
be the divisor systems
\[P_k= (p; u_1^{p^k}-u_1, \dots, u_n^{p^k}-u_n)\]   
in which the primitive form $\Delta(u_1, \dots, u_n)$
contains.\footnote{In our example there is only one, with $k=1$.} We
need only consider those whose index $k$ is not contained in one of
the other indexes $k_1, \dots, k_i$ as a divisor, because according to
the observation above, the form $\Delta(u_1, \dots, u_n)$ (if
contained in the system $P_k$) is also divisible by every system
$P_d$, whose index is a divisor of $k_i$. So we form the whole
number:
\[F(p)=(p^{k_1}-1)(p^{k_2}-1)\dots (p^{k_l}-1)\]
and look for the smallest prime $\nu$ different from $p$ which is not
contained in $F(p)$.\footnote{In our example $F(2)=(2^1-1)=1$, so
  $\nu=3$.} Then the domain $\Gamma_{\nu}$ of the $\nu$-th roots of
unity is the smallest which is a supplementary domain for
$(\mathfrak{G})$ with respect to the prime $p$.\footnote{This result
  is to be proved in the next few paragraphs. Note that in our example
  it is, as Kronecker pointed out, the field of cube roots of unity.}

Now \marginpar{[158]} we easily prove that $\Gamma_{\nu}$ is in fact a
supplementary domain for $(\mathfrak{G})$. If $p$ modulo $\nu$ belongs
to the exponent $k$, then $p$ decomposes in $\Gamma_{\nu}$ into
distinct prime factors of degree $k$. So $\Gamma_{\nu}$ is a
supplementary domain of $(\mathfrak{G})$ if and only if $\Delta(u_1,
\dots, u_n)$ does not contain the divisor system $P_k$. This claim
is only satisfied if the index $k$ is not a divisor of any of the $l$
numbers $k_1, k_2, \dots, k_l$. If this were however the case, then at
least one of the $l$ factors of the product $F(p)$, and hence also
$F(p)$ itself, would be divisible by $\nu$. Since $\nu$ is the
smallest prime not occurring in the product, the first part
of the claim\footnote{So $\Gamma_\nu$ is a supplementary domain.} is
proved. 

Further, if $\nu_1$ is a prime smaller than $\nu$ and if
$\Gamma_{\nu_1}$ is the field constituted by the $\nu_1$th roots of
unity, then $\Gamma_{\nu_1}$ can not be a supplementary domain of
$(\mathfrak{G})$. In fact, according to the previous assumption,
$\nu_1$ is a divisor of the product $F(p)$, which means $\nu_1$
is contained in at least one of the factors $(p^{k_i}-1)$. Therefore,
the exponent $k'$ of $p$ modulo $\nu_1$ is a divisor of one of the $l$
numbers $k_1, \dots, k_l$, and so the divisor system $P_k$ is a
divisor of the form $\Delta(u_1, \dots, u_n)$, or equivalently, $p$ is
an inessential divisor for $(\mathfrak{G})$ in the domain
$\Gamma_{\nu_1}$.\footnote{So no cyclotomic field corresponding to a
  smaller prime $\nu_1$ will do the job.}

Now\footnote{We have found a cyclotomic field that will serve as a
  supplementary domain. Now Hensel wants to show that we can take a
  specific subfield.} take $\Gamma_{\nu}$ to be the supplementary
domain of roots of unity we have just determined, of lowest
degree. Then it will contain as many [sub]domains
$\Gamma_{\nu}(\lambda)$ as the number of divisors of $\nu-1$. Namely,
if
\[\nu-1 = \lambda\mu\]
is a decomposition of $\nu-1$ in two factors, then under $\nu$ there
is contained a domain of degree $\lambda$, namely the one containing
the $\lambda$ periods of $\mu$ terms formed from the $\nu$-th roots of
unity. We should now investigate which of these period domains
$\Gamma_{\nu}(\lambda)$ is of lowest degree $\lambda$ and is still a
supplementary domain.\footnote{So $\Gamma_{\nu}(\lambda)$ is the
  unique cyclic subfield of degree $\lambda$ in the cyclotomic field
  of $\nu$-th roots of unity. The ``periods'' are those defined by
  Gauss in the last chapter of his \emph{Disquisitiones Arithmeticae};
  they give an explicit basis of the field $\Gamma_\nu(\lambda)$.}

To answer this question,\footnote{Hensel first states the result: we
  choose the largest divisor of $\nu$ that satisfies a divisibility
  condition.} we think of all the divisors of $\nu-1$ arranged
according to their value in descending order and denote them by
\[\mu_1, \mu_2, \mu_3, \dots, \mu_{\varrho},\]
so that $\mu_1=\nu-1$, $\mu_{\varrho}=1$, and
\[\mu_1>\mu_2>\mu_3> \dots > \mu_{\varrho}.\]
Substituting $p$ in $F(p)$ by each element of the sequence $p^{\mu_1},
p^{\mu_2}, \dots, p^{\mu_{\varrho}}$, we see first that \marginpar{[159]}
$F(p^{\mu_1})$ is divisible by $\nu$, because every one of the factors
\[p^{\mu_1k_l}-1=(p^{r-1})^{k_l}-1\]
contains this prime. But $F(p^{\mu_{\varrho}})=F(p)$ does not contain
this prime, according to the above assumption on $\nu$. Let then
$\mu$ be the first, and therefore largest, of these numbers such that 
\[F(p^{\mu})=(p^{\mu k_1}-1) \dots (p^{\mu k_l}-1)\]
is no longer divisible by $\nu$. Let $\lambda$ be the
complementary divisor of $\nu-1$ to $\mu$, that is,
\[\lambda \mu=\nu -1.\]
Then the domain $\Gamma_{\nu}(\lambda)$ containing the $\lambda$
periods with $\mu$ terms of the $\nu$th roots of unity, is the
smallest\footnote{This is the claim. The proof follows. In our
  example, of course, the only divisors of $3-1=2$ are $\mu_1=2$ and
  $\mu_2=1$ and the only subfield that works is $\Gamma_3$ itself.}
[subfield] which is still a supplementary domain for $(\mathfrak{G})$.

That this is a supplementary domain, one sees as follows: If $\kappa$
is the exponent of $p^{\mu}$ modulo $\nu$, that is, the smallest whole
number for which the difference $(p^{\kappa\mu}-1)$ is divisible by $\nu$,
then $p$ decomposes in $\Gamma_{\nu}(\lambda)$ into distinct prime
factors of degree $\kappa$. The domain $\Gamma_{\nu}(\lambda)$ is then
a supplementary domain of $(\mathfrak{G})$ if the primitive form
$\Delta(u_1, \dots, u_n)$ is not contained the module system
$P_\kappa$, so if $\kappa$ is not among the numbers $k_1, \dots,
k_l$. But if this were the case, then one of the factors
$p^{\mu k_i}-1$, and so the number $F(p^{\mu})$, would be divisible by
$\nu$, which is not the case.\footnote{So we have shown that the
  chosen $\Gamma_\nu(\lambda)$ is a supplementary domain. It remains
  to show that it is the smallest.}

Now taking a different period domain $\Gamma_{\nu}(\lambda')$ of
smaller degree, so that $\lambda'<\lambda$, and since
$\lambda'\mu'=\nu-1$, so that $\mu'>\mu$, then this cannot be a
supplementary domain for $(\mathfrak{G})$. Namely, if $\kappa'$ is the
exponent of $p^{\mu'}$ modulo $\nu$, then $p^{\mu'\kappa'} -1$ is
divisible by $\nu$, and $\kappa'$ must be a divisor of one of the
numbers $k_1, \dots, k_l$. Since $F(p^{\mu'})$ is divisible by $\nu$
for every $\mu'>\mu$, one of the factors $(p^{\mu k_i}-1)$ is a
multiple of $\nu$, which means the exponent $\kappa'$ belonging to
$p^{\mu'}$ modulo $\nu$ is contained in $k_i$, and from this it
follows that $p$ is still an inessential divisor of the discriminants
of $(\mathfrak{G})$ in the field of rationality
$\Gamma_\nu(\lambda')$. With this, we have proved the above
conjecture, and we can summarize the end result of all the last
investigations in the following elegant theorem: 
  	
\vspace{.5\baselineskip}

\marginpar{[160]} 
\noindent\hspace{\fill}\parbox{0.8\textwidth}{\hspace*{1em}Let
  $(\mathfrak{G})$ be a given field of degree $n$ and let $\Delta(u_1,
  \dots, u_n)$ be the discriminant of the fundamental equation freed
  from its numerical factor. Further let $p$ be any real prime and
  denote by \[P_{k_1}, P_{k_2}, \dots, P_{k_l}\]
  the divisor systems of the form
  \[P_k=(p; u_1^{p^k}-u_1, \dots, u_n^{p^k}-u_n)\]
  which $\Delta$ contains, chosen so that among their indexes
  $k_1, \dots, k_l$ none is a multiple of the others. 

  Let $\nu$ be the smallest prime which does not divide the integer
  \[F(p)= (p^{k_1}-1) \dots (p^{k_l}-l).\]
  Then the field $\Gamma_{\nu}$ of the $\nu$-th roots of unity is the
  smallest [cyclotomic field] which is a supplementary domain of the
  domain $(\mathfrak{G})$. Further, if $\mu$ is the largest divisor of
  $\nu -1$, for which the number
  \[F(p^{\mu})\]
  does not contain the prime $\nu$, the period field
  $\Gamma_{\nu}(\lambda)$ contained in it, defined by the $\lambda$
  periods of $\mu$ terms in the $\nu$th roots of unity, is the
  smallest [subfield] which still has this property.}

\vspace{\baselineskip}

Berlin, November 17, 1893

\clearpage


\end{document}